\documentclass{amsart}
\usepackage{amsmath,amssymb,amsfonts,amsthm,tikz-cd,hyperref,mathtools}
\usepackage[style=alphabetic]{biblatex}
\usepackage{tikz,tcolorbox}
\usetikzlibrary{shapes, positioning, arrows.meta, intersections}
\hypersetup{hidelinks}

\theoremstyle{definition}
\newtheorem{definition}{Definition}[section]
\newtheorem{situation}[definition]{Situation}

\theoremstyle{plain}
\newtheorem{lemma}[definition]{Lemma}
\newtheorem{proposition}[definition]{Proposition}
\newtheorem{theorem}[definition]{Theorem}
\newtheorem*{theorem*}{Theorem}
\newtheorem{corollary}[definition]{Corollary}
\theoremstyle{remark}
\newtheorem{remark}[definition]{Remark}
\newtheorem{example}[definition]{Example}

\title{Geometric invariant theory for graded additive groups}
\author{Yikun Qiao}
\date{\today}
\address{Yikun Qiao, Institute of Mathematics, Academy of Mathematics and Systems Science, Chinese Academy of Sciences, Beijing, 100190, China}
\email{yikunqiao@amss.ac.cn}

\begin{document}

\begin{abstract}
	We consider geometric invariant theory for \emph{graded additive groups}, groups of the form $\mathbb{G}_a^r\rtimes_w\mathbb{G}_m$ such that the $\mathbb{G}_m$-action on $\mathbb{G}_a^r$ is a scalar multiplication with weight $w\in\mathbb{N}_+$. We provide an algorithm of equivariant birational modifications, such that we can apply the geometric invariant theory of B\'erczi-Doran-Hawes-Kirwan. In particular, the geometric $\mathbb{G}_a^r$-quotient exists. This complements B\'erczi-Doran-Hawes-Kirwan, in the special case of one grading weight. 
\end{abstract}

\maketitle
\setcounter{tocdepth}{1}
\tableofcontents

\section{Introduction}
{
	Let $\Bbbk$ be an algebraically closed field of characteristic zero. Let $\mathbb{G}_a^r\cong\mathbb{A}^r$ denote the additive group of dimension $r\in\mathbb{N}_+$. For $w\in\mathbb{N}_+$, let $\mathbb{G}_m$ act on $\mathbb{G}_a^r$ such that the induced representation $\mathbb{G}_m$ on $\mathrm{Lie}(\mathbb{G}_a^r)$ is a scalar multiplication of weight $w$. Denote the semi-direct product by $\mathbb{G}_a^r\rtimes_w\mathbb{G}_m$. We will consider geometric invariant theory of $\mathbb{G}_a^r\rtimes_w\mathbb{G}_m$-actions. 

	{
		An \emph{affine cone (of schemes)} (\cite[\href{https://stacks.math.columbia.edu/tag/062P}{Section 062P}]{stacks-project}) refers to an affine morphism $\pi:C\to S$ such that $\pi_*\mathcal{O}_C$ is endowed with the structure of a quasi-coherent $\mathbb{N}$-graded $\mathcal{O}_S$-algebra $\pi_*\mathcal{O}_C=\bigoplus_{n\geq0}\mathcal{A}_n$ such that $\mathcal{A}_0=\mathcal{O}_S$. In particular $C\cong\underline{\mathrm{Spec}}_S\big(\bigoplus_{n\geq0}\mathcal{A}_n\big)$. We call $\pi:C\to S$ an \emph{affine cone of varieties} if $C,S$ are varieties. The \emph{zero section} of $\pi:C\to S$, denoted by $Z\hookrightarrow C$, refers to the closed subscheme $Z:=\mathbb{V}\big(\bigoplus_{n>0}\mathcal{A}_n\big)$. In particular the zero section $Z$ is the base $S$ embedded in $C$, i.e. the composition $Z\hookrightarrow C\to S$ is an isomorphism. 

		\begin{situation}\label{situation: G_a^r act on cone}
			Let $\pi:C\to S$ be an \emph{affine cone} of (non-empty) varieties. Let $\mathcal{A}:=\bigoplus_{n\geq0}\mathcal{A}_n$ with $\mathcal{A}_0=\mathcal{O}_S$ be the associated graded $\mathcal{O}_S$-algebra. Let $Z\hookrightarrow C$ denote the zero section. Let $\mathbb{G}_a^r$ act on $C$ such that $\pi:C\to S$ is $\mathbb{G}_a^r$-invariant. The affine cone structure is equivalent to a $\mathbb{G}_m$-action. Assume actions of $\mathbb{G}_a^r$ and $\mathbb{G}_m$ constitute an action of $\mathbb{G}_a^r\rtimes_w\mathbb{G}_m$. 
		\end{situation}
		\begin{remark}
			This situation covers the moduli problem $G\curvearrowright S_\beta$ of an unstable stratum $S_\beta\subseteq X^{\mathrm{us}}$ in the Hesselink-Kempf-Kirwan-Ness stratification for a GIT problem $G\curvearrowright X\subseteq\mathbb{P}^n$, in the special case when $\lambda_\beta:\mathbb{G}_m\to P_\beta$ acts on $\mathrm{Lie}(U_\beta)$ with one weight, where $U_\beta\subseteq P_\beta$ is the unipotent radical. Definitions of $\lambda_\beta,P_\beta,S_\beta$ are given later, or see \cite{zbMATH03880866} for details. 
		\end{remark}
	}

	\subsection{Quotients of unstable strata}
	{
		For a linear action of reductive group $G$ on a projective variety $X\subseteq\mathbb{P}^n$, the unstable locus has the Hesselink-Kempf-Kirwan-Ness (HKKN) stratification (cf. \cite{zbMATH03600245}, \cite{zbMATH03630927}, \cite{zbMATH03880866}, \cite{zbMATH03975073}). 
		\begin{align}
			X^{\mathrm{us}}=\bigsqcup_{\beta\in\mathcal{B}\setminus0}S_\beta,\quad S_\beta\cong G\times_{P_\beta}Y_\beta^{\mathrm{ss}}
		\end{align}
		where $Y_\beta^{\mathrm{ss}}\subseteq X$ are locally closed sub-varieties and $P_\beta\subseteq G$ are parabolic subgroups acting on $Y_\beta^{\mathrm{ss}}$. Moreover, there exists a one-parameter group $\lambda_\beta:\mathbb{G}_m\to P_\beta$ and a retraction to the fixed point subvariety (cf. \cite{zbMATH03427530}) 
		\begin{align}
			p_\beta:Y_\beta^{\mathrm{ss}}\to Z_\beta^{\mathrm{ss}}:=\big(Y_\beta^{\mathrm{ss}}\big)^{\lambda_\beta}. 
		\end{align}
		We can apply geometric invariant theory (GIT, \cite{zbMATH00467370}) to construct a projective good categorical quotient $X^{\mathrm{ss}}\to X/\!\!/G$, which is a projective completion of the geometric quotient $X^{\mathrm{s}}\to X^{\mathrm{s}}/G$. However, GIT yields the empty quotient when applied on $G\curvearrowright X^{\mathrm{us}}\subseteq \mathbb{P}^n$. To get a non-empty quotient of the unstable stratum $S_\beta$, it is natural to consider the quotient problem for $P_\beta\curvearrowright Y_\beta^{\mathrm{ss}}$. The parabolic group has a finite central extension $U_\beta\rtimes (T_\beta\times K_\beta)\to P_\beta$ where $U_\beta\subseteq P_\beta$ is the unipotent radical, $T_\beta$ is the central torus, and $K_\beta$ is the semisimple factor. To get a $P_\beta$-quotient, we can construct first a $U_\beta$-quotient and then a $T_\beta\times K_\beta$-quotient, or alternatively first a $U_\beta\rtimes\lambda_\beta$-quotient and then a $(T_\beta/\lambda_\beta)\times K_\beta$-quotient. We focus on the step of non-reductive quotients and consider 
		\begin{align}\label{equation: graded unipotent act on Y_beta^ss}
			U_\beta\rtimes\lambda_\beta\curvearrowright Y_\beta^{\mathrm{ss}},\quad p_\beta:Y_\beta^{\mathrm{ss}}\to Z_\beta^{\mathrm{ss}}. 
		\end{align}

		We observe the following facts. 
		\begin{itemize}
			\item The conjugation action of $\lambda_\beta:\mathbb{G}_m\to P_\beta$ on $U_\beta$ \emph{grades} the unipotent group, in the sense that all weights of the $\mathbb{G}_m$-representation on $\mathrm{Lie}(U_\beta)$ are positive. 

			\item The retraction $p_\beta:Y_\beta^{\mathrm{ss}}\to Z_\beta^{\mathrm{ss}}$ defines an affine cone. We can extend $p_\beta$ to $\widetilde{p_\beta}:Y_\beta(\mathbb{P}^n)\to Z_\beta(\mathbb{P}^n)$ such that $Z_\beta^{\mathrm{ss}}\subseteq Z_\beta(\mathbb{P}^n)$ is locally closed and $p_\beta:Y_\beta^{\mathrm{ss}}\to Z_\beta^{\mathrm{ss}}$ is the base change of $\widetilde{p_\beta}$. The closed sub-variety $Z_\beta(\mathbb{P}^n)\hookrightarrow \mathbb{P}^n$ is $\mathbb{P}(V_\beta)$ for certain $\lambda_\beta$-weight subspace $V_\beta\subseteq\Bbbk^{n+1}$, and the sub-variety $Y_\beta(\mathbb{P}^n)$ consists of $[v]\in\mathbb{P}^n$ such that $\lim_{t\to0}\lambda_\beta(t).[v]\in Z_\beta(\mathbb{P}^n)$. See \cite{zbMATH03880866} for details. We can check that $\widetilde{p_\beta}$ defines an affine cone, so does its base change $p_\beta$. 

			\item The retraction $p_\beta:Y_\beta^{\mathrm{ss}}\to Z_\beta^{\mathrm{ss}}$ is $U_\beta\rtimes\lambda_\beta$-invariant. The $\lambda_\beta$-invariance is easy to see. For $u\in U_\beta$ and $y\in Y_\beta^{\mathrm{ss}}$ we have $p_\beta(u.y)=\lim_{t\to0}\lambda_\beta(t).u.y=\big(\lim_{t\to0}\lambda_\beta(t)u\lambda_\beta(t)^{-1}\big).\big(\lim_{t\to0}\lambda_\beta(t).y\big)=\lim_{t\to0}\lambda_\beta(t).y=p_\beta(y)$. This proves the $U_\beta$-invariance. 
		\end{itemize}
		According to these facts, we need to consider actions of the form $U\rtimes\mathbb{G}_m\curvearrowright\pi:C\to S$ where: 
		\begin{itemize}
			\item $U$ is a unipotent group graded by $\mathbb{G}_m$, in the sense that all weights of $\mathbb{G}_m\curvearrowright\mathrm{Lie}(U)$ are positive; 
			\item $\pi:C\to S$ is an affine cone; 
			\item $U\rtimes\mathbb{G}_m$ acts on $C$ such that $\pi:C\to S$ is invariant, and the $\mathbb{G}_m$-action corresponds to the affine cone structure. 
		\end{itemize}
		It is easy to see that Situation \ref{situation: G_a^r act on cone} is the special case when $\mathbb{G}_m\curvearrowright\mathrm{Lie}(U)$ has only one weight. For the non-reductive quotient problem \eqref{equation: graded unipotent act on Y_beta^ss}, Situation \ref{situation: G_a^r act on cone} covers the special case when the conjugate action of $\lambda_\beta$ on $\mathrm{Lie}(U_\beta)$ has only one weight. The case with multiple weights is more complicated and will be studied in the future. 

		An important application of GIT is the moduli problem of sheaves, semistable or unstable. Simpson in \cite{zbMATH00714749} constructed moduli schemes of semistable sheaves using GIT. Also, it is noticed that sheaves of fixed Harder-Narasimhan type are bounded. In \cite{zbMATH06109728}, Hoskins and Kirwan showed that the moduli problem of sheaves of fixed Harder-Narasimhan type can be viewed as the quotient problem of a refined unstable stratum. Moduli problems of objects of fixed Harder-Narasimhan type in other linear Abelian categories are also examples of quotient problems of unstable strata, for instance \cite{zbMATH06353674}, \cite{zbMATH07020384}, \cite{hoskins2021quotients}, \cite{hamilton2022instability}. In particular, modulo difficulties in quotients by reductive groups, these moduli problems of unstable objects all correspond to quotients of non-reductive groups on an affine cone similar to Situation \ref{situation: G_a^r act on cone}, where unipotent radicals are not necessarily Abelian. 
	}

	\subsection{Current states}
	{
		Non-reductive GIT has been studied in \cite{zbMATH06963807}, \cite{bérczi2020projective}, \cite{hoskins2021quotients}, e.t.c. for \emph{graded unipotent groups} $\hat U:=U\rtimes\mathbb{G}_m$, where $U$ is unipotent and the $\mathbb{G}_m$-action on $\mathrm{Lie}(U)$ has all weights positive. Let $\hat{U}$ act on an affine cone $\pi:C\to S$ similar to Situation \ref{situation: G_a^r act on cone}. In \cite{zbMATH06963807}, geometric quotients by $U$ and $\hat{U}$ were constructed, under the condition called ``semistability coincides with stability'', that is $\dim\mathrm{Stab}_U(z)=0$ for all $z\in Z$, equivalently $\dim\mathrm{Stab}_U(x)=0$ for all $x\in C$. It was observed in \cite[Remark 7.12]{bérczi2020projective} the condition can be relaxed to that $\dim\mathrm{Stab}_{U_i}$ is constant for certain chain of subgroups $U=U_0\supseteq U_1\supseteq\cdots\supseteq U_k=\{e\}$, or further relaxed to 
		\begin{align}
			\min_{x\in C}\dim\mathrm{Stab}_{U_i}(x)=\min_{z\in Z}\dim\mathrm{Stab}_{U_i}(z),\quad i=0,1,\cdots,k
		\end{align}
		if we allow a birational equivariant modification. The case when 
		\begin{align}
			\min_{x\in C}\dim\mathrm{Stab}_{U_i}(x)<\min_{z\in Z}\dim\mathrm{Stab}_{U_i}(z)
		\end{align}
		for some $i$ were not considered previously. It was conjectured in \cite{bérczi2020projective} that we can blow up $\overline{UZ}$ iteratively to expect the relaxed ``semistability coincides with stability'' holds eventually. This paper implements this method in Situation \ref{situation: G_a^r act on cone}, by proving termination of this sequence of equivariant birational modifications. 

		In \cite{jackson2021moduli}, the GIT of \cite{bérczi2020projective} already sufficed to construct moduli schemes of sheaves of Harder-Narasimhan length 2. In \cite{zbMATH07790942}, GIT for $\hat{U}$ has been applied to hyperbolicity questions to get a polynomial degree bound for the Green-Griffiths-Lang conjecture. 
	}

	\subsection{Main result and layout}
	{
		Let $\mathbb{G}_a^r\rtimes_w\mathbb{G}_m$ act on $\pi:C\to S$ as in Situation \ref{situation: G_a^r act on cone}. The desired condition is that $\dim\mathrm{Stab}_{\mathbb{G}_a^r}$ is constant on $C$ (cf. \cite{bérczi2020projective}, \cite{hoskins2021quotients}). The main result is Theorem \ref{theorem: the algorithm is valid}, which proves that a sequence of modifications terminate to the desired situation, defining an algorithm (Figure \ref{figure: full flowchart}). The test condition $e(\pi)=\infty$ in Figure \ref{figure: full flowchart} is equivalent to the set-theoretic condition 
		\begin{align}\label{equation: dimStab attains min on Z}
			\min_{x\in C}\dim\mathrm{Stab}_{\mathbb{G}_a^r}(x)=\min_{z\in Z}\dim\mathrm{Stab}_{\mathbb{G}_a^r}(z). 
		\end{align}

		We first introduce the notion of \emph{modification}. Let $V\subsetneq C$ be a $\mathbb{G}_a^r\rtimes_w\mathbb{G}_m$-invariant closed subvariety. In this paper, the \emph{modification of $\pi:C\to S$ along $V$} refers to $\pi':C'\to S'$ in the commutative diagram 
		\begin{align}
			\begin{tikzcd}[ampersand replacement=\&]
				C'\ar[r,"p"]\ar[d,"\pi'"]\&C\ar[d,"\pi"]\\
				S'\ar[r,"q"]\&S
			\end{tikzcd}
		\end{align}
		where $C'\subseteq\mathrm{Bl}_V(C)$ is the open stratum of the Bia\l ynicki-Birula stratification (cf. \cite{zbMATH03432438}, \cite{zbMATH07128154}) for the $\mathbb{G}_m$-action on $\mathrm{Bl}_V(C)$, and $\pi':C'\to S'$ is an affine cone with a $\mathbb{G}_a^r\rtimes_w\mathbb{G}_m$-action satisfying Situation \ref{situation: G_a^r act on cone}. We call $V\subsetneq C$ the \emph{centre} of the modification. It is easy to see that $p:C'\to C$ is $\mathbb{G}_a^r\rtimes_w\mathbb{G}_m$-equivariant and birational, and that $q:S'\to S$ is projective. 

		We need two kinds of modifications, described in Section \ref{section: modification I} and Section \ref{section: modification II}. For modification of type I, the centre is 
		\begin{align}
			\overline{\mathbb{G}_a^r.(Z\setminus Z^{\mathrm{reg}})}\subsetneq C
		\end{align}
		where $Z^{\mathrm{reg}}:=\{z\in Z:\dim\mathrm{Stab}_{\mathbb{G}_a^r}(z)\textrm{ is the minimum on }Z\}$. It can happen that $Z=Z^{\mathrm{reg}}$ and the modification does nothing. For modification of type II, we only consider when $\dim\mathrm{Stab}_{\mathbb{G}_a^r}$ is constant on $Z$ and when $\min_{x\in C}\dim\mathrm{Stab}_{\mathbb{G}_a^r}(x)<\dim\mathrm{Stab}_{\mathbb{G}_a^r}(z)$. Under these conditions, we choose the centre to be 
		\begin{align}
			\mathbb{G}_a^r.Z=\overline{\mathbb{G}_a^r.Z}\subsetneq C
		\end{align}
		where the $\mathbb{G}_a^r$-sweep is closed by Lemma \ref{lemma: dimStab const on Z implies U.Z closed}. 

		The main theorem is the following. 
		{
			\begin{theorem*}[Theorem \ref{theorem: the algorithm is valid}]
				Let $\mathbb{G}_a^r\rtimes_w\mathbb{G}_m$ act on $\pi:C\to S$ as in Situation \ref{situation: G_a^r act on cone}. Then the following steps give a well-defined algorithm, terminating in finite time. 
				\begin{itemize}
					\item[(i)] Input $\pi:C\to S$. 
					\item[(ii)] Run Modification I, the modification with centre $\overline{\mathbb{G}_a^r.(Z\setminus Z^{\mathrm{reg}})}\subsetneq C$. 
					\item[(iii)] If $\dim\mathrm{Stab}_{\mathbb{G}_a^r}$ attains minimum on $Z$ (cf. \eqref{equation: dimStab attains min on Z}), then output $\pi:C\to S$. \\Else, run Modification II, the modification with centre $\mathbb{G}_a^r.Z\subsetneq C$, and then return (ii). 
				\end{itemize}
				Moreover, the output $\pi:C\to S$ with the $\mathbb{G}_a^r\rtimes_w\mathbb{G}_m$-action satisfies that $\dim\mathrm{Stab}_{\mathbb{G}_a^r}$ is constant on $C$. In particular there is a $\mathbb{G}_m$-equivariant geometric $\mathbb{G}_a^r$-quotient $C\to C/\mathbb{G}_a^r$. 
			\end{theorem*}
		}

		In Section \ref{section: prelim}, we prepare several notions. There is a coherent $\mathcal{O}_C$-module $\mathcal{Q}(\pi)$ such that Fitting loci of $\mathcal{Q}(\pi)$ stratify $C$ according to $\dim\mathrm{Stab}_{\mathbb{G}_a^r}$ (Proposition \ref{proposition: dimStab and Fit}). We also define an index $(d(\pi),e(\pi))\in\mathbb{N}\times\overline{\mathbb{N}_+}$ (Definition \ref{definition: index (d,e)}), which will bound the looping number in Figure \ref{figure: full flowchart}. We also recall the important Slice Theorem \ref{theorem: slice theorem} for the study of locally nilpotent derivations. Finally in this section, we describe an example of $\mathbb{G}_a^r$-quotient (Example \ref{example: G_a^r act on linear scheme over Grass with b-dimensional orbit}), which is stronger than a geometric quotient but weaker than a principal bundle. It proves to be a nice local model. 

		In Section \ref{section: main theorems}, we first review GIT for $\mathbb{G}_a^r\rtimes_w\mathbb{G}_m$ under the condition that $\dim\mathrm{Stab}_{\mathbb{G}_a^r}$ is constant (Theorem \ref{theorem: GIT with UU}). In this case $C\to C/\mathbb{G}_a^r$ exists, which is locally the base change of Example \ref{example: G_a^r act on linear scheme over Grass with b-dimensional orbit}. Next we describe the algorithm and prove the main result Theorem \ref{theorem: the algorithm is valid}, with details of modifications postponed to later sections. 

		In Section \ref{section: modification I} and Section \ref{section: modification II}, we study two types of modifications. After modification of type I, we attain the constancy of $\dim\mathrm{Stab}_{\mathbb{G}_a^r}$ on $Z$ (Proposition \ref{proposition: modification for constant dimStab on Z}). Moreover, both modifications do not increase the index $(d(\pi),e(\pi))$, and modification of type II will reduce it strictly (Corollary \ref{corollary: index change in modification I}, Proposition \ref{proposition: modification II}). They are the keys for the termination. Some technical arguments are separated in the appendix. 
	}

	\subsection{A working example}
	{
		Here is an example that we can go through the algorithm. We describe affine cones in each step below and verify these descriptions later in Section \ref{subsection: working example - Mod I} for Modifications of type I and in Section \ref{subsection: working example - Mod II} for Modifications of type II. We also work out the associated $\Omega_{\mathcal{A}/\mathcal{A}_0}\to\mathfrak{u}^*\otimes\mathcal{A}$ and $\mathrm{Fit}_d(\mathcal{Q}(\pi))$ in Example \ref{example: working example - Q and Fit}, and the indices $(d(\pi),e(\pi))$ in Example \ref{example: wokring example - index}. 
		\begin{example}\label{example: working example}
			Choose $\rho,\sigma\in\mathbb{N}_+$ with $0<\sigma<\rho$. Let $w\in\mathbb{N}_+$. Consider 
			\begin{equation}
				C:=\bigg\{\begin{pmatrix}\begin{pmatrix}a_{11}\\a_{21}&a_{22}\end{pmatrix},e,\begin{pmatrix}f_1\\f_2\end{pmatrix}\end{pmatrix}\in\mathbb{A}^6:a_{ij},e,f_i\in\Bbbk\bigg\}\cong\mathbb{A}^6
			\end{equation}
			and for simplicity denote $A:=\begin{pmatrix}a_{11}\\a_{21}&a_{22}\end{pmatrix}$ and $f:=\begin{pmatrix}f_1\\f_2\end{pmatrix}$. Then a point is represented by $(A,e,f)\in C$. Consider the action of $\mathbb{G}_a^2\rtimes_w\mathbb{G}_m$ 
			\begin{equation}
				\begin{split}
					u&.(A,e,f):=(A,e,f+Au),\quad u:=\begin{pmatrix}u_1\\u_2\end{pmatrix}\in\mathbb{G}_a^2\\
					t&.(A,e,f):=\begin{pmatrix}\begin{pmatrix}1\\&t^\rho\end{pmatrix}A\;,\;t^\sigma e\;,\;\begin{pmatrix}t^w\\&t^{\rho+w}\end{pmatrix}f\end{pmatrix},\quad t\in\mathbb{G}_m. 
				\end{split}
			\end{equation}
			We have an affine cone $\pi:C\to S$ with $Z\cong S\cong\{(A,e,f):a_{11}\in\Bbbk,\textrm{ others}=0\}\cong\mathbb{A}^1$. We can check 
			\begin{equation}
				\min_{x\in C}\dim\mathrm{Stab}_{\mathbb{G}_a^2}(x)=0,\quad\min_{z\in Z}\dim\mathrm{Stab}_{\mathbb{G}_a^2}(z)=1. 
			\end{equation}
			Let $\pi^{(0)}:C^{(0)}\to S^{(0)}$ denote $\pi:C\to S$. The algorithm runs as follows 
			\begin{equation}
				\begin{tikzcd}
					\pi^{(0)}\ar[d,"\textrm{I}"]&\pi^{(2)}\ar[d,equal,"\textrm{I}"]&\pi^{(4)}\ar[d,"\textrm{I}"]\\
					\pi^{(1)}\ar[ru,"\textrm{II}"]&\pi^{(3)}\ar[ru,"\textrm{II}"]&\pi^{(5)}. 
				\end{tikzcd}
			\end{equation}
			Affine cones are: 
			\begin{itemize}
				\item $C^{(0)}=C=\mathrm{Spec}\big(\Bbbk[a_{11},a_{21},a_{22},e,f_1,f_2]\big)$ and $S^{(0)}=S=\mathrm{Spec}\big(\Bbbk[a_{11}]\big)$; 

				\item $C^{(1)}=\mathrm{Spec}\big(\Bbbk\big[a_{11},\frac{a_{21}}{a_{11}},\frac{a_{22}}{a_{11}},\frac{e}{a_{11}},\frac{f_1}{a_{11}},\frac{f_2}{a_{11}}\big]\big)$ and $S^{(1)}=\mathrm{Spec}\big(\Bbbk[a_{11}]\big)$; 

				\item $C^{(3)}=C^{(2)}=\mathrm{Spec}\big(\Bbbk\big[a_{11},\frac{a_{21}}{e},\frac{a_{22}}{e},\frac{e}{a_{11}},\frac{f_1}{a_{11}},\frac{f_2}{e}\big]\big)$\\
				and $S^{(3)}=S^{(2)}=\mathrm{Spec}\big(\Bbbk[a_{11}]\big)$; 

				\item $C^{(4)}=$
				\begin{equation*}\begin{split}
					{} & \textstyle{\mathrm{Spec}\big(\Bbbk\big[a_{11},\frac{a_{21}}{e},\frac{a_{22}}{a_{21}},\frac{e^2}{a_{11}a_{21}},\frac{f_1}{a_{11}},\frac{f_2}{a_{21}}\big]\big)}\\
					\textstyle{\bigcup}{} & \textstyle{\mathrm{Spec}\big(\Bbbk\big[a_{11},\frac{a_{21}}{a_{22}},\frac{a_{22}}{e},\frac{e^2}{a_{11}a_{22}},\frac{f_1}{a_{11}},\frac{f_2}{a_{22}}\big]\big)}
				\end{split}\end{equation*}
				and $S^{(4)}=\mathrm{Spec}\big(\Bbbk\big[a_{11},\frac{a_{22}}{a_{21}}\big]\big)\bigcup \mathrm{Spec}\big(\Bbbk\big[a_{11},\frac{a_{21}}{a_{22}}\big]\big)$; 

				\item $C^{(5)}=$
				\begin{equation*}
					\begin{split}
						{} & \textstyle{\mathrm{Spec}\big(\Bbbk\big[a_{11},\frac{a_{21}^2}{a_{22}e},\frac{a_{22}}{a_{21}},\frac{e^2}{a_{11}a_{22}},\frac{f_1}{a_{11}},\frac{f_2}{a_{22}}-\frac{f_1a_{21}}{a_{11}a_{22}}\big]\big)}\\
						\textstyle{\bigcup}{} & \textstyle{\mathrm{Spec}\big(\Bbbk\big[a_{11},\frac{a_{21}}{a_{22}},\frac{a_{22}}{e},\frac{e^2}{a_{11}a_{22}},\frac{f_1}{a_{11}},\frac{f_2}{a_{22}}\big]\big)}; 
					\end{split}
				\end{equation*}
				and $S^{(5)}=\mathrm{Spec}\big(\Bbbk\big[a_{11},\frac{a_{22}}{a_{21}}\big]\big)\bigcup \mathrm{Spec}\big(\Bbbk\big[a_{11},\frac{a_{21}}{a_{22}}\big]\big)$
			\end{itemize}
			Centres are: 
			\begin{itemize}
				\item $\overline{\mathbb{G}_a^2.\big(Z^{(0)}\setminus Z^{(0),\mathrm{reg}}\big)}=\mathbb{V}\big(a_{11},a_{21},a_{22},e,f_1,f_2\big)$, i.e. $0\in\mathbb{A}^6$; 
				\item $\mathbb{G}_a^2.Z^{(1)}=\mathbb{V}\big(\frac{a_{21}}{a_{11}},\frac{a_{22}}{a_{11}},\frac{e}{a_{11}},\frac{f_2}{a_{11}}\big)$; 
				\item $\overline{\mathbb{G}_a^2.\big(Z^{(2)}\setminus Z^{(2),\mathrm{reg}}\big)}=\emptyset$; 
				\item $\mathbb{G}_a^2.Z^{(3)}=\mathbb{V}\big(\frac{a_{21}}{e},\frac{a_{22}}{e},\frac{e}{a_{11}},\frac{f_2}{e}\big)$; 
				\item $\overline{\mathbb{G}_a^2.\big(Z^{(4)}\setminus Z^{(4),\mathrm{reg}}\big)}=\mathbb{V}\big(\frac{a_{21}}{e},\frac{a_{22}}{a_{21}},\frac{e^2}{a_{11}a_{21}},\frac{f_1a_{22}}{a_{11}a_{21}},\frac{f_2}{a_{21}}-\frac{f_1}{a_{11}}\big)\bigcup\emptyset$. 
			\end{itemize}
		\end{example}
	}

	\subsection*{Convention}
	{
		A \emph{variety} refers to a reduced irreducible scheme, separated of finite type over $\Bbbk$. A \emph{point} of a variety usually refers to a closed point, or equivalently a $\Bbbk$-point. The phrase $x\in C$ refers to that $x$ is a (closed/$\Bbbk$-rational) point of $C$. We use $\mathbb{V}$ for both closed subschemes associated to an ideal sheaf (cf. Example \ref{example: working example}) or linear schemes associated to a locally free sheaf (cf. Example \ref{example: G_a^r act on linear scheme over Grass with b-dimensional orbit}). 
	}

	\subsection*{Acknowledgement}
	{
		The author is supported by the National Key R\&D Program of China (No2021YFA1002300).
	}
}

\section{Preliminaries}\label{section: prelim}
{
	Let $\mathbb{G}_a^r\times_w\mathbb{G}_m$ act on $\pi:C\to S$ as in Situation \ref{situation: G_a^r act on cone}. Let $\mathcal{A}=\bigoplus_{n\geq0}\mathcal{A}_n$ denote the associated graded $\mathcal{O}_S$-algebra, such that $\mathcal{A}_0=\mathcal{O}_S$ and $C=\underline{\mathrm{Spec}}_S(\mathcal{A})$. Let $\mathfrak{u},\mathfrak{u}^*$ denote the Lie algebra of $\mathbb{G}_a^r$ and its dual. 
	\subsection{Lie algebra actions}
	{
		In the following diagram, horizontal morphisms are immersions 
		\begin{align}
			\begin{tikzcd}[ampersand replacement=\&]
				C\ar[r,"\imath"]\ar[d,equal]\&\mathbb{G}_a^r\times C\ar[d,"\psi"]\\
				C\ar[r,"\Delta"]\&C\times_SC
			\end{tikzcd}\quad\begin{cases}
				\imath:&x\mapsto (e,x)\\
				\psi:&(u,x)\mapsto (u.x,x)\\
				\Delta:&x\mapsto (x,x). 
			\end{cases}
		\end{align}
		It induces a morphism of conormal sheaves $\Omega_{\mathcal{A}/\mathcal{A}_0}\to \mathfrak{u}^*\otimes\mathcal{A}$ by \cite[\href{https://stacks.math.columbia.edu/tag/01R4}{Tag 01R4}]{stacks-project}. Define 
		\begin{align}
			\mathcal{Q}(\pi):=\mathrm{coker}\big(\Omega_{\mathcal{A}/\mathcal{A}_0}\to \mathfrak{u}^*\otimes\mathcal{A}\big)
		\end{align}

		\begin{remark}
			When $C$ is smooth, the action $\mathbb{G}_a^r\curvearrowright C$ induces an action of $\mathfrak{u}$, represented by $\mathfrak{u}\to T_{C/S}$, where $T_{C/S}$ refers to the relative tangent bundle. In this case, the morphism $\Omega_{\mathcal{A}/\mathcal{A}_0}\to \mathfrak{u}^*\otimes\mathcal{A}$ above is dual to $\mathfrak{u}\to T_{C/S}$. Choose dual bases $\xi_1,\cdots,\xi_r\in\mathfrak{u}$ and $u_1,\cdots,u_r\in\mathfrak{u}^*$. Each $\xi_i$ induces a vector field on $C$, still denoted by $\xi_i$. The morphism $\Omega_{\mathcal{A}/\mathcal{A}_0}\to \mathfrak{u}^*\otimes\mathcal{A}$ is $\omega\mapsto \sum_{i=1}^ru_i\otimes\langle\omega,\xi_i\rangle$, where $\langle\omega,\xi_i\rangle$ is the contraction of the 1-form $\omega$ with the tangent vector field $\xi_i$.

		\end{remark}

		\begin{remark}
			In the proof the Proposition \ref{proposition: dimStab and Fit}, we will show that the fibre of the quotient $\mathfrak{u}^*\otimes\mathcal{A}\to \mathcal{Q}(\pi)$ at $x\in C$ is equivalent to the quotient $\mathfrak{u}^*\to \mathrm{Stab}_{\mathfrak{u}}(x)^*$. 
		\end{remark}

		\begin{example}\label{example: working example - Q and Fit}
			Assume $\pi^{(i)}:C^{(i)}\to S^{(i)}$ of Example \ref{example: working example} are known. We describe $\Omega_{\mathcal{A}^{(i)}/\mathcal{A}^{(i)}_0}\to\mathfrak{u}^*\otimes\mathcal{A}^{(i)}$ and $\mathrm{Fit}_d(\mathcal{Q}(\pi^{(i)}))$. Note that always $\mathrm{Fit}_{-1}(\mathcal{Q}(\pi^{(i)}))=0$ and $\mathrm{Fit}_2(\mathcal{Q}(\pi^{(i)}))=\mathcal{A}$. 
			\begin{itemize}
				\item $\mathcal{A}^{(0)}=\Bbbk[a_{11},a_{21},a_{22},e,f_1,f_2]$ and $\mathcal{A}^{(0)}_0=\Bbbk[a_{11}]$. The morphism of conormal sheaves is 
				\begin{equation*}
					\begin{split}
						\Omega_{\mathcal{A}^{(0)}/\mathcal{A}^{(0)}_0}&\to \mathfrak{u}^*\otimes\mathcal{A}^{(0)},\quad \mathrm{d}a_{21},\mathrm{d}a_{22},\mathrm{d}e\mapsto 0\\ 
						\mathrm{d}f_1&\mapsto u_1\otimes a_{11}\\
						\mathrm{d}f_2&\mapsto u_1\otimes a_{21}+u_2\otimes u_2. 
					\end{split}
				\end{equation*}
				Modulo differentials in the kernel, it is represented by the matrix $\big(\begin{smallmatrix}a_{11}\\ a_{21}&a_{22}\end{smallmatrix}\big)$ 
				\begin{equation*}
					\begin{pmatrix}
						\mathrm{d}f_1\\
						\mathrm{d}f_2
					\end{pmatrix}\mapsto \begin{pmatrix}
						a_{11}\\ 
						a_{21}&a_{22}
					\end{pmatrix}\begin{pmatrix}
						u_1\otimes 1\\ 
						u_2\otimes 1
					\end{pmatrix}. 
				\end{equation*}
				Then $\mathrm{Fit}_d(\mathcal{Q}(\pi^{(0)}))$ is generated by all $(2-d)\times(2-d)$-minors of $\big(\begin{smallmatrix}a_{11}\\ a_{21}&a_{22}\end{smallmatrix}\big)$. The Fitting ideals are 
				\begin{equation*}
					\mathrm{Fit}_0(\mathcal{Q}(\pi^{(0)}))=\langle a_{11}a_{22}\rangle,\quad\mathrm{Fit}_1(\mathcal{Q}(\pi^{(0)}))=\langle a_{11},a_{21},a_{22}\rangle. 
				\end{equation*}

				\item $\mathcal{A}^{(1)}=\Bbbk\big[a_{11},\frac{a_{21}}{a_{11}},\frac{a_{22}}{a_{11}},\frac{e}{a_{11}},\frac{f_1}{a_{11}},\frac{f_2}{a_{11}}\big]$ and $\mathcal{A}^{(1)}_0=\Bbbk[a_{11}]$. The morphism is 
				\begin{equation*}
					\begin{split}
						\Omega_{\mathcal{A}^{(1)}/\mathcal{A}^{(1)}_0}&\to \mathfrak{u}^*\otimes\mathcal{A}^{(1)},\quad \mathrm{d}\tfrac{a_{21}}{a_{11}},\mathrm{d}\tfrac{a_{22}}{a_{11}},\mathrm{d}\tfrac{e}{a_{11}}\mapsto0\\
						\begin{pmatrix}
							\mathrm{d}\tfrac{f_1}{a_{11}}\\ 
							\mathrm{d}\tfrac{f_2}{a_{11}}
						\end{pmatrix}&\mapsto \begin{pmatrix}
							1\\ 
							\tfrac{a_{21}}{a_{11}}&\tfrac{a_{22}}{a_{11}}
						\end{pmatrix}\begin{pmatrix}
							u_1\otimes 1\\ 
							u_2\otimes 1
						\end{pmatrix}. 
					\end{split}
				\end{equation*}
				The Fitting ideals are 
				\begin{equation*}
					\mathrm{Fit}_0(\mathcal{Q}(\pi^{(1)}))=\langle \tfrac{a_{22}}{a_{11}}\rangle,\quad\mathrm{Fit}_1(\mathcal{Q}(\pi^{(1)}))=\mathcal{A}^{(1)}.
				\end{equation*}

				\item $\mathcal{A}^{(3)}=\mathcal{A}^{(2)}=\Bbbk\big[a_{11},\frac{a_{21}}{e},\frac{a_{22}}{e},\frac{e}{a_{11}},\frac{f_1}{a_{11}},\frac{f_2}{e}\big]$ and $\mathcal{A}^{(3)}_0=\mathcal{A}^{(2)}_0=\Bbbk[a_{11}]$. The morphism is ($i=2,3$)
				\begin{equation*}
					\begin{split}
						\Omega_{\mathcal{A}^{(i)}/\mathcal{A}^{(i)}_0}&\to \mathfrak{u}^*\otimes\mathcal{A}^{(i)},\quad \mathrm{d}\tfrac{a_{21}}{e},\mathrm{d}\tfrac{a_{22}}{e},\mathrm{d}\tfrac{e}{a_{11}}\mapsto0\\
						\begin{pmatrix}
							\mathrm{d}\tfrac{f_1}{a_{11}}\\ 
							\mathrm{d}\tfrac{f_2}{e}
						\end{pmatrix}&\mapsto \begin{pmatrix}
							1\\ 
							\tfrac{a_{21}}{e}&\tfrac{a_{22}}{e}
						\end{pmatrix}\begin{pmatrix}
							u_1\otimes 1\\ 
							u_2\otimes 1
						\end{pmatrix}. 
					\end{split}
				\end{equation*}
				The Fitting ideals are 
				\begin{equation*}
					\mathrm{Fit}_0(\mathcal{Q}(\pi^{(i)}))=\langle \tfrac{a_{22}}{e}\rangle,\quad\mathrm{Fit}_1(\mathcal{Q}(\pi^{(i)}))=\mathcal{A}^{(i)}.
				\end{equation*}

				\item On two affine opens 
				\begin{itemize}
					\item $\mathcal{A}^{(4),1}=\Bbbk\big[a_{11},\frac{a_{21}}{e},\frac{a_{22}}{a_{21}},\frac{e^2}{a_{11}a_{21}},\frac{f_1}{a_{11}},\frac{f_2}{a_{21}}\big]$ and $\mathcal{A}^{(4),1}_0=\Bbbk\big[a_{11},\frac{a_{22}}{a_{21}}\big]$. The morphism is 
					\begin{equation*}
						\begin{split}
							\Omega_{\mathcal{A}^{(4),1}/\mathcal{A}^{(4),1}_0}&\to \mathfrak{u}^*\otimes\mathcal{A}^{(4),1},\quad \mathrm{d}\tfrac{a_{21}}{e},\mathrm{d}\tfrac{e^2}{a_{11}a_{21}}\mapsto0\\
							\begin{pmatrix}
								\mathrm{d}\tfrac{f_1}{a_{11}}\\ 
								\mathrm{d}\tfrac{f_2}{a_{21}}
							\end{pmatrix}&\mapsto \begin{pmatrix}
								1\\ 
								1&\tfrac{a_{22}}{a_{21}}
							\end{pmatrix}\begin{pmatrix}
								u_1\otimes 1\\ 
								u_2\otimes 1
							\end{pmatrix}. 
						\end{split}
					\end{equation*}
					The Fitting ideals are 
					\begin{equation*}
						\mathrm{Fit}_0(\mathcal{Q}(\pi^{(4),1}))=\langle \tfrac{a_{22}}{a_{21}}\rangle,\quad\mathrm{Fit}_1(\mathcal{Q}(\pi^{(4),1}))=\mathcal{A}^{(4),1}.
					\end{equation*}

					\item $\mathcal{A}^{(4),2}=\Bbbk\big[a_{11},\frac{a_{21}}{a_{22}},\frac{a_{22}}{e},\frac{e^2}{a_{11}a_{22}},\frac{f_1}{a_{11}},\frac{f_2}{a_{22}}\big]$ and $\mathcal{A}^{(4),2}_0=\Bbbk\big[a_{11},\frac{a_{21}}{a_{22}}\big]$. The morphism is 
					\begin{equation*}
						\begin{split}
							\Omega_{\mathcal{A}^{(4),2}/\mathcal{A}^{(4),2}_0}&\to \mathfrak{u}^*\otimes\mathcal{A}^{(4),2},\quad \mathrm{d}\tfrac{a_{22}}{e},\mathrm{d}\tfrac{e^2}{a_{11}a_{21}}\mapsto0\\
							\begin{pmatrix}
								\mathrm{d}\tfrac{f_1}{a_{11}}\\ 
								\mathrm{d}\tfrac{f_2}{a_{22}}
							\end{pmatrix}&\mapsto \begin{pmatrix}
								1\\ 
								\tfrac{a_{21}}{a_{22}}&1
							\end{pmatrix}\begin{pmatrix}
								u_1\otimes 1\\ 
								u_2\otimes 1
							\end{pmatrix}. 
						\end{split}
					\end{equation*}
					The Fitting ideals are 
					\begin{equation*}
						\mathrm{Fit}_0(\mathcal{Q}(\pi^{(4),2}))=\mathcal{A}^{(4),2},\quad\mathrm{Fit}_1(\mathcal{Q}(\pi^{(4),2}))=\mathcal{A}^{(4),2}.
					\end{equation*}
				\end{itemize}

				\item On two affine opens 
				\begin{itemize}
					\item $\mathcal{A}^{(5),1}=\Bbbk\big[a_{11},\frac{a_{21}^2}{a_{22}e},\frac{a_{22}}{a_{21}},\frac{e^2}{a_{11}a_{22}},\frac{f_1}{a_{11}},\frac{f_2}{a_{22}}-\frac{f_1a_{21}}{a_{11}a_{22}}\big]$ and $\mathcal{A}^{(5),1}_0=\Bbbk\big[a_{11},\frac{a_{22}}{a_{21}}\big]$. The morphism is 
					\begin{equation*}
						\begin{split}
							\Omega_{\mathcal{A}^{(5),1}/\mathcal{A}^{(5),1}_0}&\to \mathfrak{u}^*\otimes\mathcal{A}^{(5),1},\quad \mathrm{d}\tfrac{a_{21}^2}{a_{22}e},\mathrm{d}\tfrac{e^2}{a_{11}a_{22}}\mapsto0\\
							\begin{pmatrix}
								\mathrm{d}\tfrac{f_1}{a_{11}}\\ 
								\mathrm{d}\big(\tfrac{f_2}{a_{22}}-\tfrac{f_1a_{21}}{a_{11}a_{22}}\big)
							\end{pmatrix}&\mapsto \begin{pmatrix}
								1\\ 
								&1
							\end{pmatrix}\begin{pmatrix}
								u_1\otimes 1\\ 
								u_2\otimes 1
							\end{pmatrix}. 
						\end{split}
					\end{equation*}
					The Fitting ideals are 
					\begin{equation*}
						\mathrm{Fit}_0(\mathcal{Q}(\pi^{(5),1}))=\mathcal{A}^{(5),1},\quad\mathrm{Fit}_1(\mathcal{Q}(\pi^{(5),1}))=\mathcal{A}^{(5),1}.
					\end{equation*}

					\item $\mathcal{A}^{(5),2}=\mathcal{A}^{(4),2}$. See items for $\mathcal{A}^{(4),2}$. 
				\end{itemize}
			\end{itemize}
		\end{example}
	}

	\subsection{Unipotent stabilisers}
	{
		For $x\in C$, let $\mathrm{Stab}_{\mathfrak{u}}(x)\subseteq\mathfrak{u}$ denote the Lie algebra stabiliser. Dimensions of $\mathrm{Stab}_{\mathfrak{u}}(x)$ are determined by fibres of $\mathcal{Q}(\pi)$. 
		{
			\begin{proposition}\label{proposition: dimStab and Fit}
				Let $x\in C$ be a closed point with ideal $\mathcal{I}_x\subseteq\mathcal{A}$. Let $d\in\mathbb{N}$. The following are equivalent: 
				\begin{itemize}
					\item $\dim\mathrm{Stab}_{\mathfrak{u}}(x)>d$; 
					\item $\mathrm{Fit}_d(\mathcal{Q}(\pi))\subseteq\mathcal{I}_x$. 
				\end{itemize}
			\end{proposition}
			\begin{proof}
				The statement is local. We can assume $C$ is affine. Write $A,I_x,Q(\pi)$ for $\mathcal{A}_x,\mathcal{I}_x,\mathcal{Q}(\pi)$. Let $\kappa(x)$ denote the residue field of $x\in C$, i.e. $\kappa(x)=A/I_x$. Denote $V^\vee:=\mathrm{Hom}_{\kappa(x)}(V,\kappa(x))$ for $\kappa(x)$-spaces. The following sequence is exact 
				\begin{align}
					\begin{tikzcd}[ampersand replacement=\&]
						0\ar[r]\&\big(Q(\pi)\otimes\kappa(x)\big)^\vee\ar[r]\&\big(\mathfrak{u}^*\otimes\kappa(x)\big)^\vee\ar[r]\&\big(\Omega_{A/A_0}\otimes\kappa(x)\big)^\vee. 
					\end{tikzcd}
				\end{align}
				Let $\xi\in\mathfrak{u}$. Then $\xi\in\mathrm{Stab}_{\mathfrak{u}}(x)$ if and only if the vector field induced by $\xi$ has a zero at $x$, which is equivalent to that $\langle\omega,\xi\rangle$ has a zero at $x$ for any $\omega\in\Omega_{A/A_0}$, that is, the composition $\Omega_{A/A_0}\otimes\kappa(x)\to \mathfrak{u}^*\otimes\kappa(x)\xrightarrow{\xi\otimes\mathrm{id}_{\kappa(x)}}\kappa(x)$ is zero. The composition being zero is equivalent to 
				\begin{equation}
					\xi\otimes\mathrm{id}_{\kappa(x)}\in \ker\Big(\big(\mathfrak{u}^*\otimes\kappa(x)\big)^\vee\to \big(\Omega_{A/A_0}\otimes\kappa(x)\big)^\vee\Big)\cong \big(Q(\pi)\otimes\kappa(x)\big)^\vee. 
				\end{equation}
				This proves 
				\begin{align}
					\mathrm{Stab}_{\mathfrak{u}}(x)\cong \big(Q(\pi)\otimes\kappa(x)\big)^\vee. 
				\end{align}
				Then $\dim\mathrm{Stab}_{\mathfrak{u}}(x)>d$ if and only if $\dim Q(\pi)\otimes\kappa(x)>d$, which is equivalent to $\mathrm{Fit}_d(Q(\pi))\subseteq I_x$ by \cite[\href{https://stacks.math.columbia.edu/tag/07ZC}{Tag 07ZC}]{stacks-project}. 

				We provide another proof of the last sentence, to show how Fitting ideals naturally emerge. Choose a surjective map $A^{\oplus J}\to \Omega_{A/A_0}$. Let $M\in A^{r\times J}$ denote the matrix representing the composition $A^{\oplus J}\to \Omega_{A/A_0}\to\mathfrak{u}^*\otimes A$. For $\xi\in\mathfrak{u}$, let $\underline{\xi}\in A^{1\times r}$ denote the matrix representing the induced map $\mathfrak{u}^*\otimes A\xrightarrow{\xi\otimes\mathrm{id}_A} A$. For $M\in A^{r\times J}$ and $\underline{\xi}\in A^{1\times r}$, let $M(x)\in\kappa(x)^{r\times J}$ and $\underline{\xi}(x)\in\kappa(x)^{1\times r}$ denote the images along $A\to \kappa(x)$. Note that $\xi\mapsto\underline{\xi}(x)$ is a linear isomorphism. With these notations, we have 
				\begin{equation}
					\xi\in\mathrm{Stab}_{\mathfrak{u}}(x)\quad \iff \quad \underline{\xi}(x)M(x)=0. 
				\end{equation}

				If $\dim\mathrm{Stab}_{\mathfrak{u}}(x)>d$, then there exist linearly independent $\xi_0,\cdots,\xi_d\in\mathrm{Stab}_{\mathfrak{u}}(x)$. Then $\underline{\xi_i}(x)M(x)=0$, which implies $\mathrm{rank}(M(x))<r-d$. Then all $(r-d)\times(r-d)$-minors of $M(x)$ vanish, which implies that all $(r-d)\times (r-d)$-minors of $M$ are in $I_x\subseteq A$, that is $\mathrm{Fit}_d(Q(\pi))\subseteq I_x$ by definition. 

				If $\mathrm{Fit}_d(Q(\pi))\subseteq I_x$, then all $(r-d)\times(r-d)$-minors of $M$ are in $I_x$. Then all $(r-d)\times(r-d)$-minors of $M(x)$ are zero, which is equivalent to $\mathrm{rank}(M(x))<r-d$. Then there exist linearly independent $y_0,\cdots,y_d\in\kappa(x)^{1\times r}$ such that $y_iM(x)=0$. Since $\mathfrak{u}\to \kappa(x)^{1\times r}$ is an isomorphism, we can find linearly independent $\xi_0,\cdots,\xi_d\in \mathfrak{u}$ such that $\underline{\xi_i}(x)=y_i$. Then $\xi_i\in\mathrm{Stab}_{\mathfrak{u}}(x)$, which implies $\dim\mathrm{Stab}_{\mathfrak{u}}(x)>d$. 
			\end{proof}
		}

		{
			Let $\overline{\mathbb{N}_+}$ denote the totally ordered set $\mathbb{N}_+\sqcup\{\infty\}$. Consider the lexicographic order on $\mathbb{N}\times\overline{\mathbb{N}_+}$, i.e. $(d,e)<(d',e')$ if and only if $d<d'$ or $\begin{cases}d=d'\\e<e'\end{cases}$. 
			\begin{definition}\label{definition: index (d,e)}
				We define an index for the $\mathbb{G}_a^r\rtimes_w\mathbb{G}_m$-action on $\pi:C\to S$ in Situation \ref{situation: G_a^r act on cone}
				\begin{align}
					(d(\pi),e(\pi))\in\mathbb{N}\times\overline{\mathbb{N}_+}
				\end{align}
				via 
				\begin{align}
					d(\pi)&:=\min_{z\in Z}\dim\mathrm{Stab}_{\mathfrak{u}}(z)\\
					e(\pi)&:=\begin{cases}
						\min\{n\in\mathbb{N}:\mathrm{Fit}_{d(\pi)-1}(\mathcal{Q}(\pi))\cap \mathcal{A}_n\ne0\},&\textrm{if }\mathrm{Fit}_{d(\pi)-1}(\mathcal{Q}(\pi))\ne0\\
						\infty,&\textrm{if }\mathrm{Fit}_{d(\pi)-1}(\mathcal{Q}(\pi))=0. 
					\end{cases}
				\end{align}
			\end{definition}
			Note that: 
			\begin{itemize}
				\item $d(\pi)$ is the integer such that 
				\begin{equation}
					\begin{cases}
						\dim\mathrm{Stab}_{\mathfrak{u}}(z)>d(\pi)-1\textrm{ for all }z\in Z\\ 
						\dim\mathrm{Stab}_{\mathfrak{u}}(z)\leq d(\pi)\textrm{ for some }z\in Z. 
					\end{cases}
				\end{equation}
				By Proposition \ref{proposition: dimStab and Fit} and $\mathcal{A}_{>0}=\bigcap_{z\in Z}\mathfrak{m}_z$, we have 
				\begin{equation}\label{equation: Fit definition of d(pi)}
					\mathrm{Fit}_{d(\pi)-1}(\mathcal{Q}(\pi))\subseteq\mathcal{A}_{>0}\quad\textrm{and}\quad\mathrm{Fit}_{d(\pi)}(\mathcal{Q}(\pi))\not\subseteq\mathcal{A}_{>0}. 
				\end{equation}
				\item We can not have $e(\pi)=0$, since $\mathrm{Fit}_{d(\pi)-1}(\mathcal{Q}(\pi))\subseteq\mathcal{A}_{>0}$ and $e(\pi)$ is the minimal degree of elements of $\mathrm{Fit}_{d(\pi)-1}(\mathcal{Q}(\pi))$. 

				\item $e(\pi)=\infty$ if and only if $d(\pi)=\min_{x\in C}\dim\mathrm{Stab}_{\mathfrak{u}}(x)$. By definition $e(\pi)=\infty$ if and only if $\mathrm{Fit}_{d(\pi)-1}(\mathcal{Q}(\pi))=0$. Then $\dim\mathrm{Stab}_{\mathfrak{u}}(x)\geq d(\pi)$ for all $x\in C$ by Proposition \ref{proposition: dimStab and Fit}. 
			\end{itemize}

			\begin{remark}\label{remark: three cases}
				There are three non-exclusive cases: 
				\begin{itemize}
					\item[(i)] $\dim\mathrm{Stab}_{\mathfrak{u}}$ are constant, i.e. $\dim\mathrm{Stab}_{\mathfrak{u}}(z)=d(\pi)$ for all $z\in Z$ and $e(\pi)=\infty$. 
					\item[(ii)] Minimal $\dim\mathrm{Stab}_{\mathfrak{u}}$ are achieved on $Z$, i.e. $e(\pi)=\infty$. 
					\item[(iii)] The opposite of (ii) holds, i.e. $e(\pi)<\infty$. 
				\end{itemize}
				Case (i) is the best case when we can get a nice $\mathbb{G}_a^r$-quotient directly (See Section \ref{section: GIT with Upstairs Unipotent Stabiliser condition} and Theorem \ref{theorem: GIT with UU}). Case (ii) is the relaxed ``semistability coincides with stability'' considered in \cite{bérczi2020projective}, \cite{hoskins2021quotients}, when one equivariant birational modification leads to case (i). Case (iii) were not considered previously. We will prove that a sequence of modifications transform case (iii) to case (i). 
			\end{remark}

			\begin{example}\label{example: wokring example - index}
				We calculate $(d(\pi),e(\pi))$ for $\pi:C\to S$ in Example \ref{example: working example}. From Example \ref{example: working example - Q and Fit}, we have $\mathrm{Fit}_0(\mathcal{Q}(\pi))=\langle a_{11}a_{22}\rangle$ and $\mathrm{Fit}_1(\mathcal{Q}(\pi))=\langle a_{11},a_{21},a_{22}\rangle$. In particular $\mathrm{Fit}_0(\mathcal{Q}(\pi))\subseteq\mathcal{A}_{>0}$ since $a_{11}a_{22}\in\mathcal{A}_\rho$, and $\mathrm{Fit}_1(\mathcal{Q}(\pi))\not\subseteq\mathcal{A}_{>0}$ since $a_{11}\in\mathcal{A}_0$. By the Fitting-ideal definition of $d(\pi)$ in \eqref{equation: Fit definition of d(pi)}, we have $d(\pi)=1$. The index $e(\pi)$ is the minimal degree of $\mathrm{Fit}_0(\mathcal{Q}(\pi))$, which is the degree of $a_{11}a_{22}$. Then $e(\pi)=\rho$. In Example \ref{example: working example - Q and Fit}, the Fitting ideals were given and we similarly get 
				\begin{itemize}
					\item $(d(\pi^{(0)}),e(\pi^{(0)}))=(1,\rho)$; 
					\item $(d(\pi^{(1)}),e(\pi^{(1)}))=(1,\rho)$; 
					\item $(d(\pi^{(3)}),e(\pi^{(3)}))=(d(\pi^{(2)}),e(\pi^{(2)}))=(1,\rho-\sigma)$; 
					\item $(d(\pi^{(4)}),e(\pi^{(4)}))=(0,\infty)$; 
					\item $(d(\pi^{(5)}),e(\pi^{(5)}))=(0,\infty)$. 
				\end{itemize}
			\end{example}

		}
	}

	\subsection{The slice theorem}
	{

		For a commutative $\Bbbk$-algebra $R$, let $\mathrm{LND}(R)$ denote the set of all \emph{locally nilpotent derivations}, that is, $\Bbbk$-linear maps $\partial:R\to R$ satisfying the Leibniz rule $\partial(fg)=f\partial(g)+\partial(f)g$. For $\partial,\partial'\in\mathrm{LND}(R)$, we have $[\partial,\partial']:=\partial\partial'-\partial'\partial\in\mathrm{LND}(R)$. We say $\partial_1,\cdots,\partial_m\in\mathrm{LND}(R)$ \emph{commute} if $[\partial_i,\partial_j]=0$ for all $1\leq i,j\leq m$. 

		Choose dual bases $\xi_1,\cdots,\xi_r\in\mathfrak{u}$ and $u_1,\cdots,u_r\in\mathfrak{u}^*$. Let $\mathbb{G}_a^r$ act on $\mathrm{Spec}(R)$ with co-action map $\sigma:R\to R[u_1,\cdots,u_r]$. Each $\xi_i\in\mathfrak{u}$ induces $\xi_i\in\mathrm{LND}(R)$ via 
		\begin{align}
			\begin{tikzcd}[ampersand replacement=\&]
				R\ar[r,"\sigma"]\&R[u_1,\cdots,u_r]\ar[r,"\partial/\partial u_i"]\&R[u_1,\cdots,u_r]\ar[rr,"u_1=\cdots=u_r=0"]\&\&R. 
			\end{tikzcd}
		\end{align}
		It is easy to see $\xi_1,\cdots,\xi_r\in\mathrm{LND}(R)$ commute. Conversely for commuting $\xi_1,\cdots,\xi_r\in\mathrm{LND}(R)$, we can define a $\mathbb{G}_a^r$-action via 
		\begin{align}
			R\to R[u_1,\cdots,u_r],\quad f\mapsto \sum_{n\in\mathbb{N}^r}\frac{\xi_1^{n_1}\cdots\xi_r^{n_r}(f)}{n_1!\cdots n_r!}u_1^{n_1}\cdots u_r^{n_r}. 
		\end{align}
		A $\mathbb{G}_a^r$-action on $\mathrm{Spec}(R)$ is equivalent to commuting $\xi_1,\cdots,\xi_r\in\mathrm{LND}(R)$. In particular, we have 
		\begin{align}
			R^{\mathbb{G}_a^r}=R^{\xi_1,\cdots,\xi_r}:=\{f\in R:\xi_1(f)=\cdots=\xi_r(f)=0\}. 
		\end{align}

		{
			The slice theorem (cf. \cite[Corollary 1.26]{zbMATH06755655}) plays an important role in the study of locally nilpotent derivations. We state a suitable version of the slice theorem. 
			\begin{theorem}[Slice Theorem]\label{theorem: slice theorem}
				Let $R$ be a commutative $\Bbbk$-algebra. Let $\xi_1,\cdots,\xi_r\in\mathrm{LND}(R)$ commute. Assume that there exist $f_1,\cdots,f_r\in R$ such that $\xi_i(f_j)=\delta_{i,j}$ for $1\leq i,j\leq r$. Then 
				\begin{itemize}
					\item $f_1,\cdots,f_r$ are algebraically independent over $R^{\xi_1,\cdots,\xi_r}$; 
					\item $R=R^{\xi_1,\cdots,\xi_r}[f_1,\cdots,f_r]$. 
				\end{itemize}
			\end{theorem}
		}
	}



	\subsection{An example of \texorpdfstring{$\mathbb{G}_a^r$}{G\_a\^r}-quotient}
	{
		Consider the Grassmannian of $b$-dimensional quotients of $\mathfrak{u}$ and the total space of the universal quotient over it, denoted by $\pi:\mathbb{V}(\mathcal{U}^\vee)\to\mathrm{Grass}(\mathfrak{u},b)$. We will describe a $\mathbb{G}_a^r$-action on the total space $\mathbb{V}(\mathcal{U}^\vee)$ such that $\mathrm{Grass}(\mathfrak{u},b)$ is the $\mathbb{G}_a^r$-quotient. This example provides a local model for $\mathbb{G}_a^r\curvearrowright C$ when $\dim\mathrm{Stab}_{\mathfrak{u}}=r-b$ is constant. 

		When $\dim\mathrm{Stab}_{\mathfrak{u}}=r-b$, the $\mathbb{G}_a^r$-quotient was constructed in \cite{zbMATH06963807}. We will rephrase their theorem and provide another proof in Theorem \ref{theorem: GIT with UU}. We observe that the $\mathbb{G}_a^r$-quotient $\varphi:C\to C/\mathbb{G}_a^r$ is the base change of $\mathbb{V}(\mathcal{U}^\vee)\to\mathrm{Grass}(\mathfrak{u},b)$, up to a Zariski covering of $C/\mathbb{G}_a^r$. See Theorem \ref{theorem: GIT with UU} and its proof for details. 

		Let $0\leq b\leq r$. For $u\in\Bbbk^r$, $f\in\Bbbk^b$ and $A\in\Bbbk^{b\times r}$, the assignment $(u,f,A)\mapsto (f+Au,A)$ defines a $\mathbb{G}_a^r$-action on $\mathbb{A}^b\times\mathbb{A}^{b\times r}$. The group $\mathrm{GL}(b)$ acts naturally on $\mathbb{A}^b\times\mathbb{A}^{b\times r}$, such that the $\mathrm{GL}(b)$-action commutes with the $\mathbb{G}_a^r$-action. Consider the open subset where $\mathrm{rank}(A)=b$. We have a diagram such that vertical morphisms are geometric $\mathbb{G}_a^r$-quotients 
		\begin{align}
			\begin{tikzcd}[ampersand replacement=\&]
				\mathbb{A}^b\times\mathbb{A}^{b\times r}_{\textrm{full rank}}\ar[r]\ar[d]\&\big(\mathbb{A}^b\times\mathbb{A}^{b\times r}_{\textrm{full rank}}\big)\big/\mathrm{GL}(b)\ar[d]\\
				\mathbb{A}^{b\times r}_{\textrm{full rank}}\ar[r]\&\mathbb{A}^{b\times r}_{\textrm{full rank}}\big/\mathrm{GL}(b). 
			\end{tikzcd}
		\end{align}
		The right morphism is a linear scheme over a Grassmannian, described as follows. 
		{
			\begin{example}\label{example: G_a^r act on linear scheme over Grass with b-dimensional orbit}
				We describe a $\mathbb{G}_a^r$-action on a variety with $b$-dimensional orbits. Let $d:=r-b$ be the dimension of stabilisers. Choose dual bases 
				\begin{align}
					\xi_1,\cdots,\xi_r\in\mathfrak{u},\quad u_1,\cdots,u_r\in\mathfrak{u}^*. 
				\end{align}
				Consider $\mathrm{Gr}:=\mathrm{Grass}(\mathfrak{u},b)$, the Grassmannian of $b$-dimensional quotients of $\mathfrak{u}$. The universal quotient fits in the canonical sequence 
				\begin{align}
					\begin{tikzcd}[ampersand replacement=\&]
						0\ar[r]\&\mathcal{K}\ar[r,"i"]\&\mathfrak{u}\otimes\mathcal{O}_{\mathrm{Gr}}\ar[r,"q"]\&\mathcal{U}\ar[r]\&0
					\end{tikzcd}
				\end{align}
				with $\mathrm{rank}(\mathcal{K})=d$ and $\mathrm{rank}(\mathcal{U})=b$. The linear scheme associated to $\mathcal{U}^\vee:=\mathcal{H}\!\mathit{om}_{\mathcal{O}_{\mathrm{Gr}}}(\mathcal{U},\mathcal{O}_{\mathrm{Gr}})$ (cf. \cite[\href{https://stacks.math.columbia.edu/tag/01M1}{Tag 01M1}]{stacks-project}) is 
				\begin{align}
					\pi:\mathbb{V}(\mathcal{U}^\vee)\to\mathrm{Gr},\quad \mathbb{V}(\mathcal{U}^\vee):=\underline{\mathrm{Spec}}_{\mathrm{Gr}}\big(\mathrm{Sym}^\bullet\mathcal{U}^\vee\big). 
				\end{align}
				There is a $\mathbb{G}_a^r$-action on $\mathbb{V}(\mathcal{U}^\vee)$ such that $\pi:\mathbb{V}(\mathcal{U}^\vee)\to\mathrm{Gr}$ is invariant. The action morphism $\mathbb{G}_a^r\times\mathbb{V}(\mathcal{U}^\vee)\to\mathbb{V}(\mathcal{U}^\vee)$ corresponds to the morphism of $\mathcal{O}_{\mathrm{Gr}}$-modules 
				\begin{align}
					\phi:\mathcal{U}^\vee\to \big(\mathrm{Sym}^\bullet\mathcal{U}^\vee\big)[u_1,\cdots,u_r]:=\bigoplus_{n\geq0}\big(\mathrm{Sym}^n\mathcal{U}^\vee\big)[u_1,\cdots,u_r]
				\end{align}
				such that $\phi=\phi_0+\phi_1+\cdots$ for 
				\begin{itemize}
					\item $\phi_0$ the composition $\mathcal{U}^\vee\xrightarrow{q^\vee}\mathfrak{u}^*\otimes\mathcal{O}_{\mathrm{Gr}}\subseteq\mathcal{O}_{\mathrm{Gr}}[u_1,\cdots,u_r]$; 
					\item $\phi_1:\mathcal{U}^\vee\to \mathcal{U}^\vee[u_1,\cdots,u_r]$ the natural inclusion; 
					\item $\phi_n=0$ for $n\geq 2$. 
				\end{itemize}
				Moreover, there is an affine cone structure on $\pi:\mathbb{V}(\mathcal{U}^\vee)\to\mathrm{Gr}$ and Situation \ref{situation: G_a^r act on cone} is satisfied. We have 
				\begin{align}
					\mathcal{Q}(\pi)\cong\mathcal{K}^\vee\otimes_{\mathcal{O}_{\mathrm{Gr}}}\mathrm{Sym}^\bullet\mathcal{U}^\vee
				\end{align}
				which is locally free of rank $d$ on $\mathbb{V}(\mathcal{U}^\vee)$. Then $\dim\mathrm{Stab}_{\mathfrak{u}}(x)=d$ for all $x\in \mathbb{V}(\mathcal{U}^\vee)$ by \cite[\href{https://stacks.math.columbia.edu/tag/0C3G}{Tag 0C3G}]{stacks-project} and Proposition \ref{proposition: dimStab and Fit}. 

				We write the action in local coordinates. Let $I\subseteq\{1,\cdots,r\}$ be a subset of size $b$, with elements $i_1<\cdots<i_b$. Let $\mathrm{Gr}_I\subseteq\mathrm{Gr}$ denote the open subset of $[\mathfrak{u}\to Q]\in\mathrm{Gr}$ such that $\bigoplus_{i\in I}\Bbbk\xi_i\to \mathfrak{u}\to Q$ is an isomorphism. Then $\mathrm{Gr}_I\cong\mathbb{A}^{bd}$ with coordinates $a_{i,j}\in\mathcal{O}_{\mathrm{Gr}_I}$ for $i\in\{1,\cdots,b\}$ and $j\in\{1,\cdots,r\}\setminus I$. We set $a_{i,i_j}=\delta_{i,j}$ and write the coordinates as a $b\times r$-matrix $A:=(a_{ij})$. The matrix $A$ represents the universal quotient on $\mathrm{Gr}_I$ 
				\begin{align}
					\begin{tikzcd}[ampersand replacement=\&]
						\mathfrak{u}\otimes\mathcal{O}_{\mathrm{Gr}_I}\ar[r,"q|_{\mathrm{Gr}_I}"]\ar[d,"\cong"]\&\mathcal{U}|_{\mathrm{Gr}_I}\ar[r]\ar[d,"\cong"]\&0\\
						(\mathcal{O}_{\mathrm{Gr}_I})^{\oplus r}\ar[r,"A"]\&(\mathcal{O}_{\mathrm{Gr}_I})^{\oplus b}\ar[r]\&0. 
					\end{tikzcd}
				\end{align}
				For $1\leq i\leq b$, let $f_i:\mathcal{U}|_{\mathrm{Gr}_I}\to\mathcal{O}_{\mathrm{Gr}_I}$ denote the composition $\mathcal{U}|_{\mathrm{Gr}_I}\cong (\mathcal{O}_{\mathrm{Gr}_I})^{\oplus b}\xrightarrow{\mathrm{pr}_i}\mathcal{O}_{\mathrm{Gr}_I}$. Then $f_i\in \Gamma(\mathrm{Gr}_I,\mathcal{U}^\vee)$ together with $a_{i,j}$ form a set of coordinates of $\mathbb{V}(\mathcal{U}^\vee)_I:=\pi^{-1}(\mathrm{Gr}_I)$. We have $\phi(f_i)=\phi_0(f_i)+\phi_1(f_i)$ for 
				\begin{align}
					\phi_0(f_i)=\sum_{j=1}^ra_{ij}u_i\in\mathcal{O}_{\mathrm{Gr}}[u_1,\cdots,u_r],\quad\phi_1(f_i)=f_i\in\mathcal{U}^\vee[u_1,\cdots,u_r]. 
				\end{align}

				A point of $\mathbb{V}(\mathcal{U}^\vee)_I$ can be represented as 
				\begin{align}
					(f,A)\in \mathbb{V}(\mathcal{U}^\vee)_I,\quad f:=\begin{pmatrix}
						f_1\\\vdots\\f_b
					\end{pmatrix},\quad A:=\begin{pmatrix}
						a_{1,1}&\cdots&a_{1,r}\\
						\vdots&\cdots&\vdots\\
						a_{b,1}&\cdots&a_{b,r}
					\end{pmatrix}
					\textrm{such that }a_{i,i_j}=\delta_{i,j}. 
				\end{align}
				The action of $\mathbb{G}_a^r$ on $\mathbb{V}(\mathcal{U}^\vee)_I$ in this representation is 
				\begin{align}
					u.(f,A)=(f+Au,A),\quad u:=\begin{pmatrix}
						u_1\\\vdots\\u_r
					\end{pmatrix}\in\mathbb{G}_a^r. 
				\end{align}
			\end{example}
		}

		\begin{remark}
			The morphism $\mathbb{V}(\mathcal{U}^\vee)_I\to\mathrm{Gr}_I$ for $I=\{i_1<\cdots<i_b\}$ in Example \ref{example: G_a^r act on linear scheme over Grass with b-dimensional orbit} also illustrates the Slice Theorem \ref{theorem: slice theorem}. Denote $R:=\mathcal{O}\big(\mathbb{V}(\mathcal{U}^\vee)_I\big)$. There are functions $a_{i,j},f_i\in R$ as in Example \ref{example: G_a^r act on linear scheme over Grass with b-dimensional orbit}. If we choose dual bases $\xi_1,\cdots,\xi_r\in\mathfrak{u}$ and $u_1,\cdots,u_r\in\mathfrak{u}^*$, then we have $\xi_{i_j}(f_i)=a_{i,i_j}=\delta_{i,j}$ for $1\leq i,j\leq b$, the condition of the Slice Theorem. In particular, we have $\mathcal{O}(\mathrm{Gr}_I)=R^{\xi_{i_1},\cdots,\xi_{i_b}}$ and $f_1,\cdots,f_b$ are algebraically independent over $\mathcal{O}(\mathrm{Gr}_I)$. The conclusion can also be checked directly. 
		\end{remark}

		{
			\begin{remark}\label{remark: G_a^r act on linear scheme over Grass with b-dimensional orbit}
				There is a natural isomorphism $\mathrm{Grass}(\mathfrak{u}^*,d)\cong\mathrm{Grass}(\mathfrak{u},b)$. The canonical sequence on $\mathrm{Grass}(\mathfrak{u}^*,d)$ is the dual of that on $\mathrm{Grass}(\mathfrak{u},b)$. We reserve $\mathcal{U}$ for the universal quotient on $\mathrm{Grass}(\mathfrak{u},b)$. We tend to use $\mathrm{Grass}(\mathfrak{u}^*,d)$ in the rest of the paper, though equivalent. 
			\end{remark}
		}
	}
}

\section{Main theorems}\label{section: main theorems}
{
	We first review a theorem of GIT for $\mathbb{G}_a^r\rtimes_w\mathbb{G}_m$ in \cite{bérczi2020projective}, when case (i) of Remark \ref{remark: three cases} holds. In this case we obtain a geometric $\mathbb{G}_a^r$-quotient. Next we describe an algorithm for general cases, involving a sequence of modifications of type I or II, to achieve case (i). 

	\subsection{GIT with Upstairs Unipotent Stabiliser condition}\label{section: GIT with Upstairs Unipotent Stabiliser condition}
	{
		The terminology \emph{Upstairs Unipotent Stabiliser Assumption} comes from \cite[Assumption 4.40]{hoskins2021quotients}, which generalised ``semistability coincides with stability'' in \cite{zbMATH06963807}. In Situation \ref{situation: G_a^r act on cone}, this condition is case (i) of Remark \ref{remark: three cases}. Geometric invariant theory has been studied with this condition in \cite{bérczi2020projective}, \cite{jackson2021moduli}, \cite{hoskins2021quotients}. We mainly adapt the course of \cite[Section 7]{bérczi2020projective} for our situation. 

		Let $\mathbb{G}_a^r\rtimes_w\mathbb{G}_m$ act on $\pi:C\to S$ as in Situation \ref{situation: G_a^r act on cone}. In this section we assume $\dim\mathrm{Stab}_{\mathfrak{u}}(x)=d(\pi)$ for all $x\in C$, which by Proposition \ref{proposition: dimStab and Fit} is equivalent to 
		\begin{align}\label{equation: Fit for UU}
			\mathrm{Fit}_{d(\pi)-1}(\mathcal{Q}(\pi))=0,\quad \mathrm{Fit}_{d(\pi)}(\mathcal{Q}(\pi))=\mathcal{A}. 
		\end{align}

		There is a morphism 
		\begin{align}
			C\to \mathrm{Grass}(\mathfrak{u}^*,d(\pi)),\quad x\mapsto\big[\mathfrak{u}^*\to \mathrm{Stab}_{\mathfrak{u}}(x)^*\big]
		\end{align}
		which is the classifying morphism of the quotient $\begin{tikzcd}\mathfrak{u}^*\otimes\mathcal{A}\ar[r]& \mathcal{Q}(\pi)\ar[r]&0\end{tikzcd}$, and $\mathcal{Q}(\pi)$ is locally free of rank $d(\pi)$ by \cite[\href{https://stacks.math.columbia.edu/tag/0C3G}{Tag 0C3G}]{stacks-project} and \eqref{equation: Fit for UU}. This morphism factors through $S$ 
		\begin{align}
			\begin{tikzcd}[ampersand replacement=\&]
				C\ar[r]\ar[d,"\pi"]\&\mathrm{Grass}(\mathfrak{u}^*,d(\pi)). \\
				S\ar[ru]
			\end{tikzcd}
		\end{align}

		{
			\begin{theorem}[B\'erczi-Doran-Hawes-Kirwan, \cite{bérczi2020projective}]\label{theorem: GIT with UU}
				Let $\mathbb{G}_a^r\rtimes_w\mathbb{G}_m$ act on $\pi:C\to S$ as in Situation \ref{situation: G_a^r act on cone} with $\dim\mathrm{Stab}_{\mathfrak{u}}(x)=d(\pi)$ for all $x\in C$. Then there exists an affine cone of varieties $\varpi:C/\mathbb{G}_a^r\to S$ and a morphism of affine cones $\varphi:C\to C/\mathbb{G}_a^r$ 
				\begin{align}
					\begin{tikzcd}[ampersand replacement=\&]
						C\ar[r,dashed]\ar[d,"\varphi"]\ar[dd,swap,xshift=-10,bend right=20,"\pi"]\&\mathbb{V}(\mathcal{U}^\vee)\ar[d]\\
						C/\mathbb{G}_a^r\ar[r]\ar[d,"\varpi"]\&\mathrm{Grass}(\mathfrak{u}^*,d(\pi))\\
						S\ar[ru]
					\end{tikzcd}
				\end{align}
				such that $\varphi$ is locally a base change of $\mathbb{V}(\mathcal{U}^\vee)\to \mathrm{Grass}(\mathfrak{u}^*,d(\pi))$ in Example \ref{example: G_a^r act on linear scheme over Grass with b-dimensional orbit}. In particular $\varphi$ is a universal geometric $\mathbb{G}_a^r$-quotient, $\mathbb{G}_m$-equivariant, affine and smooth. 
			\end{theorem}
			\begin{proof}
				Focus locally on $S$ and we assume without loss of generality the following. 
				\begin{itemize}
					\item The base variety $S$ is affine. 
					\item There exists $\mathfrak{u}=\mathfrak{s}\oplus\mathfrak{v}$ such that $\dim\mathfrak{s}=d(\pi)$ and $\mathfrak{s}^*\otimes\mathcal{A}\cong\mathcal{Q}(\pi)$ in the diagram 
					\begin{align}
						\begin{tikzcd}[ampersand replacement=\&]
							\&\mathfrak{s}^*\otimes\mathcal{A}\ar[d]\ar[rd,"\cong"]\\
							\Omega_{\mathcal{A}/\mathcal{A}_0}\ar[r]\&\mathfrak{u}^*\otimes\mathcal{A}\ar[r]\&\mathcal{Q}(\pi)\ar[r]\&0. 
						\end{tikzcd}
					\end{align}
					\item Consequently $\Omega_{\mathcal{A}/\mathcal{A}_0}\to \mathfrak{v}^*\otimes\mathcal{A}$ is surjective. 
				\end{itemize}

				Let $b:=r-d(\pi)=\dim\mathfrak{v}$. Choose a basis $\beta_1,\cdots,\beta_b\in\mathfrak{v}$. Since $\Omega_{\mathcal{A}/\mathcal{A}_0}\to\mathfrak{v}^*\otimes\mathcal{A}$ is surjective, there exist $f_1,\cdots,f_b\in\mathcal{A}_w$ such that $\beta_i.f_j=\delta_{i,j}$ for $1\leq i,j\leq b$. By Slice Theorem \ref{theorem: slice theorem}, we have a polynomial ring extension $\mathcal{A}^{\mathfrak{v}}\subseteq\mathcal{A}$. We then prove that $\mathcal{A}^{\mathfrak{u}}\subseteq\mathcal{A}^{\mathfrak{v}}$ is an identity. Let $f\in \mathcal{A}^{\mathfrak{v}}$. Then $\mathrm{d}f\in\Omega_{\mathcal{A}/\mathcal{A}_0}$ maps to $0\in\mathfrak{v}^*\otimes\mathcal{A}$. Then the image of $\mathrm{d}f$ in $\mathfrak{u}^*\otimes\mathcal{A}$ is in the sub-module $\mathfrak{s}^*\otimes\mathcal{A}$. Mapping along the isomorphism $\mathfrak{s}^*\otimes\mathcal{A}\cong\mathcal{Q}(\pi)$, this image must be zero, i.e. $\mathrm{d}f\mapsto 0\in\mathfrak{u}^*\otimes\mathcal{A}$. Then $f\in\mathcal{A}^{\mathfrak{u}}$, the identity proved. 

				The morphism $\varphi:C\to C/\mathbb{G}_a^r$ corresponds to $\mathcal{A}^{\mathfrak{u}}\subseteq\mathcal{A}$. Then $\varphi$ is $\mathbb{G}_a^r$-invariant, $\mathbb{G}_m$-equivariant, affine and smooth. We then prove that $\varphi:C\to C/\mathbb{G}_a^r$ is a base change of $\mathbb{V}(\mathcal{U}^\vee)\to \mathrm{Grass}(\mathfrak{u}^*,d(\pi))$. We have $\mathcal{A}=\mathcal{A}^{\mathfrak{u}}[f_1,\cdots,f_b]$ by Theorem \ref{theorem: slice theorem}, and then $\Omega_{\mathcal{A}/\mathcal{A}^{\mathfrak{u}}}\cong \bigoplus_{i=1}^b\mathcal{A}\mathrm{d}f_i$ by \cite[\href{https://stacks.math.columbia.edu/tag/00RX}{Tag 00RX}]{stacks-project}. The pullback of $\mathcal{U}^\vee$ along $C\to\mathrm{Grass}(\mathfrak{u}^*,d(\pi))$ is 
				\begin{equation}\begin{split}
					\mathcal{U}^\vee|_C&\cong\ker\big(\mathfrak{u}^*\otimes\mathcal{A}\to\mathcal{Q}(\pi)\big)\\
					&\cong \mathrm{im}\big(\Omega_{\mathcal{A}/\mathcal{A}_0}\to\mathfrak{u}^*\otimes\mathcal{A}\big)\\
					&\cong \mathrm{im}\big(\Omega_{\mathcal{A}/\mathcal{A}^{\mathfrak{u}}}\to\mathfrak{u}^*\otimes\mathcal{A}\big)\\
					&\cong \mathrm{im}\Big(\bigoplus_{i=1}^b\mathcal{A}\mathrm{d}f_i\to\mathfrak{u}^*\otimes\mathcal{A}\Big). 
				\end{split}\end{equation}
				A morphism $C\to\mathbb{V}(\mathcal{U}^\vee)$ is equivalent to an $\mathcal{A}$-linear map $\mathcal{U}^\vee|_C\to \mathcal{A}$. Consider the map 
				\begin{align}
					\bigoplus_{i=1}^b\mathcal{A}\mathrm{d}f_i\to\mathcal{A},\quad \mathrm{d}f_i\mapsto f_i
				\end{align}
				which factors through $\mathcal{U}^\vee|_C$, defining a map $\mathcal{U}^\vee|_C\to \mathcal{A}$, thus a morphism $C\to \mathbb{V}(\mathcal{U}^\vee)$. Then the diagram is a base change square 
				\begin{align}
					\begin{tikzcd}[ampersand replacement=\&]
						C\ar[r]\ar[d,"\varphi"]\ar[rd,phantom,very near start,"\lrcorner"]\&\mathbb{V}(\mathcal{U}^\vee)\ar[d]\\
						C/\mathbb{G}_a^r\ar[r]\&\mathrm{Grass}(\mathfrak{u}^*,d(\pi)). 
					\end{tikzcd}
				\end{align}
				Note that the argument is valid locally. Globally, the quotient morphism $\varphi:C\to C/\mathbb{G}_a^r$ exists, corresponding to $\mathcal{A}^{\mathfrak{u}}\to \mathcal{A}$. However, there may not exist global sections of $\varphi:C\to C/\mathbb{G}_a^r$, so $\varphi$ is in general not a base change of $\mathbb{V}(\mathcal{U}^\vee)\to \mathrm{Grass}(\mathfrak{u}^*,d(\pi))$. 
			\end{proof}
		}
	}

	\subsection{Algorithm of modifications}
	{
		Let $\mathbb{G}_a^r\rtimes_w\mathbb{G}_m$ act on $\pi:C\to S$ as in Situation \ref{situation: G_a^r act on cone}. We will prove the following algorithm (See Figure \ref{figure: full flowchart}) terminates and the output satisfies the condition for Theorem \ref{theorem: GIT with UU} (UU in Figure \ref{figure: full flowchart} refers to upstairs-unipotent-stabiliser condition). 
		\begin{itemize}
			\item[(i)] Input $\pi:C\to S$. 
			\item[(ii)] Run Modification I (Section \ref{section: modification I}). 
			\item[(iii)] If $e(\pi)=\infty$, output $\pi:C\to S$. \\Else, run Modification II (Section \ref{section: modification II}) and return (ii). 
		\end{itemize}
		\begin{figure}
			\begin{tikzpicture}
				\newcommand\hmin{1 cm}
				\newcommand\wmin{3 cm}

				\node[draw, 
				trapezium, 
				trapezium left angle = 75, trapezium right angle = 105, 
				trapezium stretches, 
				minimum width = 0.75*\wmin, minimum height = 0.75*\hmin]
				(input)
				{$\pi:C\to S$};

				\node[below = 0.75of input, inner sep=0]
				(between input and mod I)
				{};

				\node[draw, 
				rectangle, 
				minimum width = \wmin, minimum height = \hmin, 
				below = 1.5 of input]
				(modification I)
				{Modification I}; 
				\draw (modification I.north west) -- ++ (-0.05*\wmin,0) |- (modification I.south west);
				\draw (modification I.north east) -- ++ (0.05*\wmin,0) |- (modification I.south east); 

				\node[draw, 
				rectangle, 
				minimum width = \wmin, minimum height = \hmin, 
				right = of modification I]
				(modification II)
				{Modification II}; 
				\draw (modification II.north west) -- ++ (-0.05*\wmin,0) |- (modification II.south west);
				\draw (modification II.north east) -- ++ (0.05*\wmin,0) |- (modification II.south east); 

				\node[draw, 
				diamond, aspect = 2.5, 
				minimum width = 0.75*\wmin, minimum height = 0.75*\hmin, 
				below = 0.75of modification I]
				(e pi=infty)
				{$e(\pi)=\infty$}; 

				\node[draw, 
				trapezium, 
				trapezium left angle = 75, trapezium right angle = 105, 
				trapezium stretches, 
				minimum width = 0.75*\wmin, minimum height = 0.75*\hmin, 
				below=1.5 of {e pi=infty}]
				(output)
				{UU};

				\draw[-latex] (input.south) -- (modification I.north);
				\draw[-latex] (modification I.south) -- (e pi=infty.north);
				\draw[-latex] (e pi=infty.south) -- (output.north) node[midway, fill=white]{Yes}; 
				\draw[-latex] (e pi=infty.east) -| (modification II.south) node[midway, fill=white]{No}; 
				\draw[-latex] (modification II.north) |- (between input and mod I); 
			\end{tikzpicture}
			\caption{Flowchart of algorithm}
			\label{figure: full flowchart}
		\end{figure}
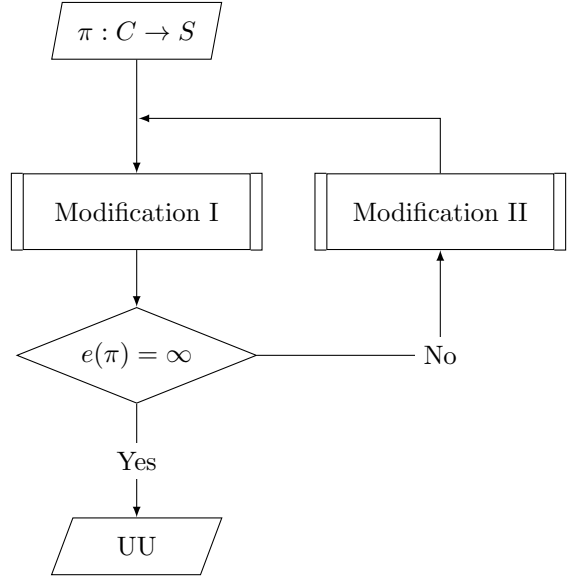

		{
			\begin{theorem}\label{theorem: the algorithm is valid}
				The algorithm above is well defined, it terminates, and the output $\pi:C\to S$ satisfies the condition of Theorem \ref{theorem: GIT with UU} 
				\begin{align}
					\dim\mathrm{Stab}_{\mathfrak{u}}(x)=d(\pi),\quad\textrm{for all }x\in C. 
				\end{align}
			\end{theorem}
			\begin{proof}
				There is no prerequisite for Modification I, so the concatenation of step (i) and step (ii) and the return from step (iii) to step (ii) are valid. To perform Modification II, we need $e(\pi)<\infty$ and $\dim\mathrm{Stab}_{\mathfrak{u}}(z)=d(\pi)$ for all $z\in Z$ (cf. Proposition \ref{proposition: modification II}). The condition $e(\pi)<\infty$ is satisfied by the logic of step (iii). By Proposition \ref{proposition: modification for constant dimStab on Z}, we have $\dim\mathrm{Stab}_{\mathfrak{u}}(z)=d(\pi)$ for all $z\in Z$ after Modification I of step (ii). Then the concatenation of step (ii) and step (iii) is valid. The algorithm is well defined. 

				Suppose the algorithm does not terminate, i.e. the loop in Figure \ref{figure: full flowchart} runs indefinitely. Modification I does not increase the index $(d(\pi),e(\pi))$ by Corollary \ref{corollary: index change in modification I}, and Modification II decreases the index $(d(\pi),e(\pi))$ strictly by Proposition \ref{proposition: modification II}. Since $\mathbb{N}\times\overline{\mathbb{N}_+}$ is well ordered, the index $(d(\pi),e(\pi))$ can not decrease indefinitely. The algorithm terminates. 

				The output $\pi:C\to S$ satisfies $e(\pi)=\infty$. Moreover, it is a result of Modification I. Then $\dim\mathrm{Stab}_{\mathfrak{u}}(z)=d(\pi)$ for all $z\in Z$ by Proposition \ref{proposition: modification for constant dimStab on Z}. Then $\pi:C\to S$ is in case (i) of Remark \ref{remark: three cases}, and equivalently $\dim\mathrm{Stab}_{\mathfrak{u}}(x)=d(\pi)$ for all $x\in C$. 
			\end{proof}

			\begin{remark}
				Let $\pi:C\to S$ be the input and $\pi':C'\to S'$ be the output. Then $d(\pi')=\min_{x\in C}\dim\mathrm{Stab}_{\mathfrak{u}}(x)$. 
			\end{remark}
		}

	}
}

\section{Modification I}\label{section: modification I}
{
	Let $\mathbb{G}_a^r\rtimes_w\mathbb{G}_m$ act on $\pi:C\to S$ as in Situation \ref{situation: G_a^r act on cone}. In this section we describe a modification 
	\begin{align}
		\begin{tikzcd}[ampersand replacement=\&]
			C'\ar[r,"p"]\ar[d,"\pi'"]\&C\ar[d,"\pi"]\\
			S'\ar[r,"q"]\&S
		\end{tikzcd}
	\end{align}
	such that $\dim\mathrm{Stab}_{\mathfrak{u}}(z')=d(\pi')=d(\pi)$ for all $z'\in Z'$, and if $\dim\mathrm{Stab}_{\mathfrak{u}}(z)=d(\pi)$ for all $z\in Z$, then $\pi=\pi'$. 

	\subsection{Centre of the modification}
	{
		Denote $Z^{\mathrm{reg}}\subseteq Z$ the open subset with minimal $\mathbb{G}_a^r$-stabilisers 
		\begin{align}
			Z^{\mathrm{reg}}:=\{z\in Z:\dim\mathrm{Stab}_{\mathfrak{u}}(z)=d(\pi)\}. 
		\end{align}
		The open immersion $Z^{\mathrm{reg}}\subseteq Z$ is the non-vanishing locus of the ideal 
		\begin{align}
			\big(\mathrm{Fit}_{d(\pi)}(\mathcal{Q}(\pi))+\mathcal{A}_{>0}\big)\big/\mathcal{A}_{>0}\subseteq\mathcal{A}/\mathcal{A}_{>0}. 
		\end{align}
		The closed subscheme $Z\setminus Z^{\mathrm{reg}}\hookrightarrow C$ is cut out by the ideal 
		\begin{align}
			\mathcal{I}:=\mathrm{Fit}_{d(\pi)}(\mathcal{Q}(\pi))+\mathcal{A}_{>0}\subseteq\mathcal{A}. 
		\end{align}

		We choose the centre of the modification as the smallest $\mathbb{G}_a^r$-invariant closed subscheme containing $Z\setminus Z^{\mathrm{reg}}$, which is the scheme theoretic image 
		\begin{align}
			\mathbb{G}_a^r\times(Z\setminus Z^{\mathrm{reg}})\hookrightarrow \mathbb{G}_a^r\times C\to C. 
		\end{align}
		The centre is thus cut out by the ideal 
		\begin{equation}
			\mathcal{J}:=\ker\big(\mathcal{A}\to \mathcal{O}(\mathbb{G}_a^r)\otimes \mathcal{A}/\mathcal{I}\big). 
		\end{equation}
	}

	\subsection{Statement and proof}
	{

		{
			\begin{proposition}\label{proposition: modification for constant dimStab on Z}
				Let $\mathbb{G}_a^r\rtimes_w\mathbb{G}_m$ act on $\pi:C\to S$ as in Situation \ref{situation: G_a^r act on cone}. Let $\pi':C'\to S'$ be the modification with centre $\overline{\mathbb{G}_a^r.(Z\setminus Z^{\mathrm{reg}})}$. Then $\dim\mathrm{Stab}_{\mathfrak{u}}(z')=d(\pi')=d(\pi)$ for all $z'\in Z'$. 
			\end{proposition}
			\begin{proof}
				We prove that $\dim\mathrm{Stab}_{\mathfrak{u}}(z')=d(\pi)$ for all $z'\in Z'$. Then by definition $d(\pi')=d(\pi)$. Let $b:=r-d(\pi)$. Without loss of generality, we assume that $S$ is affine. We have 
				\begin{align}
					C'=\bigcup_{a\in\mathcal{J}\cap \mathcal{A}_0\setminus\{0\}}C'_a,\quad C'_a:=\mathrm{Spec}\Big(\mathcal{A}\Big[\frac{\mathcal{J}}{a}\Big]\Big). 
				\end{align}
				It is easy to see that $C'$ is covered by $C'_a$ for $a\in\mathcal{J}\cap\mathcal{A}_0\setminus\{0\}$ of the following form 
				\begin{align}
					a=\det\begin{pmatrix}
						\zeta_1.f_1&\cdots&\zeta_1.f_b\\
						\vdots&\ddots&\vdots\\
						\zeta_b.f_1&\cdots&\zeta_b.f_b
					\end{pmatrix}\quad\textrm{for some }\zeta_1,\cdots,\zeta_b\in\mathfrak{u},\;f_1,\cdots,f_b\in\mathcal{A}_w. 
				\end{align}

				Consider such $C'_a$. For $i=1,\cdots,b$, denote 
				\begin{align}
					g_i:=\det\begin{pmatrix}
						\zeta_1.f_1&\cdots&\zeta_1.f_b\\
						\vdots&&\vdots\\
						f_1&\cdots&f_b\\
						\vdots&&\vdots\\
						\zeta_b.f_1&\cdots&\zeta_b.f_b
					\end{pmatrix}\in\mathcal{A}_w\quad\begin{matrix}
						\textrm{where }\begin{pmatrix}
							f_1&\cdots&f_b
						\end{pmatrix}\\
						\textrm{occupies the $i$th row. }
					\end{matrix}
				\end{align}
				We have $g_i\in\mathcal{I}$ and $\gamma_1\cdots\gamma_m.g_i\in\mathcal{I}$ for all $m\in\mathbb{N}_+$ and all $\gamma_1,\cdots,\gamma_m\in\mathfrak{u}$. Then $g_i\in\ker\big(\mathcal{A}\to \mathcal{O}(\mathbb{G}_a^r)\otimes\mathcal{A}/\mathcal{I}\big)\subseteq\mathcal{J}$ by Proposition \ref{proposition: ideal of scheme theoretic image}. Then $\frac{g_i}{a}\in\mathcal{A}\big[\frac{\mathcal{J}}{a}\big]=\mathcal{O}(C'_a)$, and then 
				\begin{align}
					1=\det\begin{pmatrix}
						\zeta_1.\frac{g_1}{a}&\cdots&\zeta_1.\frac{g_b}{a}\\
						\vdots&\ddots&\vdots\\
						\zeta_b.\frac{g_1}{a}&\cdots&\zeta_b.\frac{g_b}{a}
					\end{pmatrix}\in\mathrm{Fit}_{r-b}(\mathcal{Q}(\pi'_a))=\mathrm{Fit}_{d(\pi)}(\mathcal{Q}(\pi'_a))
				\end{align}
				which implies by Proposition \ref{proposition: dimStab and Fit}
				\begin{align}\label{equation: dimStab(x') <= d(pi)}
					\dim\mathrm{Stab}_{\mathfrak{u}}(x')\leq d(\pi),\quad\textrm{for all }x'\in C'_a. 
				\end{align}

				For all $m\in\mathbb{N}_+$ and all $\gamma_1,\cdots,\gamma_{b+1}\in\mathfrak{u}$ and all $h_1,\cdots,h_{b+1}\in\mathcal{J}^m$ we have 
				\begin{equation}\begin{split}\label{equation: Fit_(d-1) in modification I}
					&\det\begin{pmatrix}
						\gamma_1.\frac{h_1}{a^m}&\cdots&\gamma_1.\frac{h_{b+1}}{a^m}\\
						\vdots&\ddots&\vdots\\
						\gamma_{b+1}.\frac{h_1}{a^m}&\cdots&\gamma_{b+1}.\frac{h_{b+1}}{a^m}
					\end{pmatrix}\\
					=&\frac{1}{a^{m(b+1)}}\det\begin{pmatrix}
						\gamma_1.h_1&\cdots&\gamma_1.h_{b+1}\\
						\vdots&\ddots&\vdots\\
						\gamma_{b+1}.h_1&\cdots&\gamma_{b+1}.h_{b+1}
					\end{pmatrix}\in\frac{\mathcal{J}^{m(b+1)}\cap\mathrm{Fit}_{d(\pi)-1}(\mathcal{Q}(\pi))}{a^{m(b+1)}}. 
				\end{split}\end{equation}
				Determinants above generate $\mathrm{Fit}_{d(\pi)-1}(\mathcal{Q}(\pi'_a))$. Then $\mathrm{Fit}_{d(\pi)-1}(\mathcal{Q}(\pi'_a))\subseteq\mathcal{O}(C'_a)_{>0}$ since $\mathrm{Fit}_{d(\pi)-1}(\mathcal{Q}(\pi))\subseteq\mathcal{A}_{>0}$, which implies by Proposition \ref{proposition: dimStab and Fit}
				\begin{align}\label{equation: dimStab(z') > d(pi)-1}
					\dim\mathrm{Stab}_{\mathfrak{u}}(z')>d(\pi)-1,\quad\textrm{for all }z'\in Z'_a
				\end{align}
				where $Z'_a\hookrightarrow C'_a$ is cut out by $\mathcal{O}(C'_a)_{>0}$. 

				Then $\dim\mathrm{Stab}_{\mathfrak{u}}(z')=d(\pi)$ for all $z'\in Z'_a$ by \eqref{equation: dimStab(x') <= d(pi)} and \eqref{equation: dimStab(z') > d(pi)-1}. Since $C$ is covered by $C'_a$, we have $\dim\mathrm{Stab}_{\mathfrak{u}}(z')=d(\pi)$ for all $z\in Z'$. 
			\end{proof}

			The index $(d(\pi),e(\pi))$ does not increase along the modification by the following corollary. 
			\begin{corollary}\label{corollary: index change in modification I}
				Let $\pi:C\to S$ and $\pi':C'\to S'$ be as in Proposition \ref{proposition: modification for constant dimStab on Z}. Then 
				\begin{align}
					(d(\pi'),e(\pi'))\leq (d(\pi),e(\pi))\quad\textrm{in }\mathbb{N}\times\overline{\mathbb{N}_+}. 
				\end{align}
			\end{corollary}
			\begin{proof}
				We have $d(\pi)=d(\pi')$ by Proposition \ref{proposition: modification for constant dimStab on Z}. It suffices to prove $e(\pi')\leq e(\pi)$ when $e(\pi)<\infty$. 

				Consider $e(\pi)<\infty$. Since $S$ is irreducible, the index $(d(\pi),e(\pi))$ does not change if we restrict to a non-empty open subvariety. We may assume all four conditions of Lemma \ref{lemma: Fit_(d-1) cap A_e is generated by Fit_(d-1) cap Fit_(r-1) cap A_e}. For consistency, write $J=\mathcal{J}$ for the ideal of the centre. There exists $\zeta\in\mathfrak{u}$ and $h\in A_{e(\pi)+w}$ such that $\beta_1.h=\cdots=\beta_b.h=0$ and $0\ne\zeta.h\in\mathrm{Fit}_{d(\pi)-1}(\mathcal{Q}(\pi))\cap A_{e(\pi)}$ by Lemma \ref{lemma: Fit_(d-1) cap A_e is generated by Fit_(d-1) cap Fit_(r-1) cap A_e} (iii). Let $a\in J\cap A_0\setminus\{0\}$. Then $C'_a:=\mathrm{Spec}\big(A\big[\frac{J}{a}\big]\big)$ is a non-empty open subvariety of $C'$. We have 
				\begin{align}\label{equation: nonzero in Fit_(d-1)}
					0\ne\frac{\zeta.h}{a}=\det\begin{pmatrix}
						\beta_1.f_1&\cdots&\beta_1.f_b&\beta_1.\frac{h}{a}\\
						\vdots&\ddots&\vdots&\vdots\\
						\beta_b.f_1&\cdots&\beta_b.f_b&\beta_b.\frac{h}{a}\\
						\zeta.f_1&\cdots&\zeta.f_b&\zeta.\frac{h}{a}
					\end{pmatrix}\in\mathrm{Fit}_{d(\pi)-1}(\mathcal{Q}(\pi'_a))=\mathrm{Fit}_{d(\pi'_a)-1}(\mathcal{Q}(\pi'_a))
				\end{align}
				and $\frac{\zeta.h}{a}\in A_{e(\pi)}$. This proves $e(\pi')\leq e(\pi)$. 
			\end{proof}

			\begin{remark}
				The inequality of Corollary \ref{corollary: index change in modification I} is actually an equality. The inverse inequality is a consequence of \eqref{equation: Fit_(d-1) in modification I}. 
			\end{remark}
		}
	}

	\subsection{Modifications in Example \ref{example: working example}}\label{subsection: working example - Mod I}
	{
		Assume $C^{(0)}, C^{(2)}, C^{(4)}$ of Example \ref{example: working example} are known. We determine centres and describe $C^{(1)}, C^{(3)}, C^{(5)}$. We apply previous calculations from Example \ref{example: working example - Q and Fit} and Example \ref{example: wokring example - index}. 
		
		\subsubsection{$C^{(0)}$}
		{
			Geometrically, we have $Z^{(0)}=\big\{\big(\big(\begin{smallmatrix}a_{11}&0\\0&0\end{smallmatrix}\big),0,\big(\begin{smallmatrix}0\\ 0\end{smallmatrix}\big)\big)\in C^{(0)}\big\}$. Recall the $\mathbb{G}_a^2$-action is $u.(A,e,f)=(A,e,f+Au)$. Then $Z^{(0),\mathrm{reg}}=\{(A,0,0)\in Z^{(0)}:a_{11}\ne0\}$. The centre is 
			\begin{equation*}\begin{split}
				\overline{\mathbb{G}_a^2.(Z^{(0)}\setminus Z^{(0),\mathrm{reg}})}={}&\overline{\mathbb{G}_a^2.\big\{(0,0,0)\big\}}\\ 
				={}&\big\{(0,0,0)\big\}\\ 
				={}&\mathbb{V}\big(a_{11},a_{21},a_{22},e,f_1,f_2\big). 
			\end{split}\end{equation*}

			Algebraically, we have $d(\pi^{(0)})=1$ and $\mathrm{Fit}_1(\mathcal{Q}(\pi^{(0)}))=\langle a_{11},a_{21},a_{22}\rangle$. Then 
			\begin{equation*}
				\mathcal{I}^{(0)}=\mathrm{Fit}_1(\mathcal{Q}(\pi^{(0)}))+\mathcal{A}^{(0)}_{>0}=\langle a_{11}, a_{21},a_{22},e,f_1,f_2\rangle
			\end{equation*}
			and the centre is associated to the ideal 
			\begin{equation*}
				\mathcal{J}^{(0)}=\ker\big(\mathcal{A}^{(0)}\to\mathcal{O}(\mathbb{G}_a^2)\otimes \mathcal{A}^{(0)}/\mathcal{I}^{(0)}\big)=\langle a_{11}, a_{21},a_{22},e,f_1,f_2\rangle. 
			\end{equation*}
			The blowing up of $C^{(0)}$ along $\mathcal{J}^{(0)}$ is 
			\begin{equation*}
				\mathrm{Bl}_{\mathcal{J}^{(0)}}C^{(0)}=\mathrm{Proj}\big(\mathcal{A}^{(0)}\oplus\mathcal{J}^{(0)}\oplus(\mathcal{J}^{(0)})^2\oplus\cdots\big)
			\end{equation*}
			which is covered by affine opens of non-vanishing loci of $a_{11},a_{21},a_{22},e,f_1,f_2$ viewed as $\mathcal{O}(1)$-sections respectively. The open Bia\l ynicki-Birula stratum is then covered by the non-vanishing loci of sections of the minimal degree (always 0 for Modification of type I). In this case the open stratum is $\{a_{11}\ne0\}$, i.e. 
			\begin{equation*}
				C^{(1)}=\mathrm{Spec}\big(\Bbbk\big[a_{11},\tfrac{a_{21}}{a_{11}},\tfrac{a_{22}}{a_{11}},\tfrac{e}{a_{11}},\tfrac{f_1}{a_{11}},\tfrac{f_2}{a_{11}}\big]\big). 
			\end{equation*}
		}

		\subsubsection{$C^{(2)}$}
		{
			For $C^{(2)}$, we have $d(\pi^{(2)})=1$ and $\mathrm{Fit}_1(\mathcal{Q}(\pi^{(2)}))=\mathcal{A}^{(2)}$. Then $Z^{(2)}\setminus Z^{(2),\mathrm{reg}}=\emptyset$, since
			\begin{equation*}
				\mathcal{I}^{(2)}=\mathrm{Fit}_1(\mathcal{Q}(\pi^{(2)}))+\mathcal{A}^{(2)}_{>0}=\mathcal{A}^{(2)}
			\end{equation*}
			and the centre $\overline{\mathbb{G}_a^2.(Z^{(2)}\setminus Z^{(2),\mathrm{reg}})}=\emptyset$, i.e. 
			\begin{equation*}
				\mathcal{J}^{(2)}=\ker\big(\mathcal{A}^{(0)}\to\mathcal{O}(\mathbb{G}_a^2)\otimes \mathcal{A}^{(2)}/\mathcal{I}^{(2)}\big)=\mathcal{A}^{(2)}. 
			\end{equation*}
			The blowing up and taking the open Bia\l ynicki-Birula stratum are both identity morphisms 
			\begin{equation*}
				C^{(3)}=\mathrm{Bl}_{\mathcal{J}^{(2)}}C^{(2)}=C^{(2)}. 
			\end{equation*}
		}

		\subsubsection{$C^{(4)}$}
		{
			For $C^{(4)}$ we work on two of its affine opens separately. We have $d(\pi^{(4)})=0$. On $\mathrm{Spec}(\mathcal{A}^{(4),1})=\mathrm{Spec}\big(\Bbbk\big[a_{11},\tfrac{a_{21}}{e},\tfrac{a_{22}}{a_{21}},\tfrac{e^2}{a_{11}a_{21}},\tfrac{f_1}{a_{11}},\tfrac{f_2}{a_{21}}\big]\big)$, we have $\mathrm{Fit}_0(\mathcal{Q}(\pi^{(4),1}))=\langle \tfrac{a_{22}}{a_{21}}\rangle$. Then 
			\begin{equation*}
				\mathcal{I}^{(4),1}=\big\langle \tfrac{a_{22}}{a_{21}},\tfrac{a_{21}}{e},\tfrac{e^2}{a_{11}a_{21}},\tfrac{f_1}{a_{11}},\tfrac{f_2}{a_{21}}\big\rangle
			\end{equation*}
			and 
			\begin{equation*}
				\mathcal{J}^{(4),1}=\ker\big(\mathcal{A}^{(4),1}\to \mathcal{O}(\mathbb{G}_a^2)\otimes\mathcal{A}^{(4),1}/\mathcal{I}^{(4),1}\big)=\big\langle\tfrac{a_{21}}{e},\tfrac{a_{22}}{a_{21}},\tfrac{e^2}{a_{11}a_{21}},\tfrac{a_{22}f_1}{a_{11}a_{21}},\tfrac{f_2}{a_{21}}-\tfrac{f_1}{a_{11}}\big\rangle. 
			\end{equation*}
			The centre is 
			\begin{equation*}
				\overline{\mathbb{G}_a^2.(Z^{(4),1}\setminus Z^{(4),1,\mathrm{reg}})}=\mathbb{V}\big(\tfrac{a_{21}}{e},\tfrac{a_{22}}{a_{21}},\tfrac{e^2}{a_{11}a_{21}},\tfrac{a_{22}f_1}{a_{11}a_{21}},\tfrac{f_2}{a_{21}}-\tfrac{f_1}{a_{11}}\big)
			\end{equation*}
			Among five generators of $\mathcal{J}^{(4),1}$, the only one of degree 0 is $\tfrac{a_{22}}{a_{21}}$. Then 
			\begin{equation*}\begin{split}
				\mathcal{A}^{(5),1}={}&\mathcal{A}^{(4),1}\Big[\frac{\mathcal{J}^{(4),1}}{\tfrac{a_{22}}{a_{21}}}\Big]\\ 
				={}&\textstyle{\Bbbk\big[a_{11},\frac{a_{21}^2}{a_{22}e},\frac{a_{22}}{a_{21}},\frac{e^2}{a_{11}a_{22}},\frac{f_1}{a_{11}},\frac{f_2}{a_{22}}-\frac{f_1a_{21}}{a_{11}a_{22}}\big]}. 
			\end{split}\end{equation*}

			On $\mathrm{Spec}(\mathcal{A}^{(4),2})$, we have $\mathrm{Fit}_0(\mathcal{Q}(\pi^{(4),2}))=\mathcal{A}^{(4),2}$. The centre is empty by the same reason for $C^{(2)}$. Then $\mathcal{A}^{(5),2}=\mathcal{A}^{(4),2}$. 

		}
	}
}

\section{Modification II}\label{section: modification II}
{
	Let $\mathbb{G}_a^r\times_w\mathbb{G}_m$ act on $\pi:C\to S$ as in Situation \ref{situation: G_a^r act on cone}. Assume further that 
	\begin{align}\label{equation: dim Stab constant on Z}
		\dim\mathrm{Stab}_{\mathfrak{u}}(z)=d(\pi),\quad\textrm{for all }z\in Z
	\end{align}
	which is achieved by the modification in Section \ref{section: modification I}. If $e(\pi)=\infty$, then the upstairs-unipotent-stabiliser condition in Section \ref{section: GIT with Upstairs Unipotent Stabiliser condition} holds. In this section we describe a modification 
	\begin{align}
		\begin{tikzcd}[ampersand replacement=\&]
			C'\ar[r,"p"]\ar[d,"\pi'"]\&C\ar[d,"\pi"]\\
			S'\ar[r,"q"]\&S
		\end{tikzcd}
	\end{align}
	such that 
	\begin{align}
		(d(\pi'),e(\pi'))<(d(\pi),e(\pi))\quad\textrm{if}\quad e(\pi)<\infty. 
	\end{align}
	Note that the modification may not be the identity when $e(\pi)=\infty$. 

	\subsection{Centre of the modification}
	{
		We choose the centre as the smallest $\mathbb{G}_a^r$-invariant closed subscheme containing $Z$, which is the scheme theoretic image 
		\begin{align}
			\mathbb{G}_a^r\times Z\hookrightarrow\mathbb{G}_a^r\times C\to C. 
		\end{align}

		{
			\begin{lemma}\label{lemma: dimStab const on Z implies U.Z closed}
				When \eqref{equation: dim Stab constant on Z} holds, the image of $\mathbb{G}_a^r\times Z\to C$ is closed. 
			\end{lemma}
			\begin{proof}
				The statement is local on $S$, and we may assume: 
				\begin{itemize}
					\item $S$ is affine such that $\pi:C\to S$ corresponds to $A_0\to A=\bigoplus_{n\geq0}A_n$; 
					\item there exist $\beta_1,\cdots,\beta_b\in\mathfrak{u}$ and $f_1,\cdots,f_b\in A_w$ for $b:=r-d(\pi)$ such that $\beta_i.f_j=\delta_{i,j}$. 
				\end{itemize}
				Let $\mathfrak{v}:=\mathrm{Span}_\Bbbk\{\beta_1,\cdots,\beta_b\}$. Let $U'\subseteq \mathbb{G}_a^r$ be the subgroup generated by $\mathfrak{v}$. We have $\dim\mathrm{Stab}_{\mathfrak{v}}(x)=0$ for all $x\in C$, by surjectivity of $\Omega_{A/A_0}\to \mathfrak{v}^*\otimes A$ and Proposition \ref{proposition: dimStab and Fit}. We can apply Theorem \ref{theorem: GIT with UU} for the $U'\rtimes_w\mathbb{G}_m$-action. In particular $C\to C/U'$ is a geometric quotient. Then $U'.Z\subseteq C$ is closed, since it is the preimage of the closed subset $Z\subseteq C/U'$ along $C\to C/U'$. Then $U'.Z$ is the support of the scheme theoretic image of $U'\times Z\to C$ by \cite[\href{https://stacks.math.columbia.edu/tag/01R8}{Tag 01R8}]{stacks-project}. 

				It suffices to prove 
				\begin{align}
					\ker\big(A\to\mathcal{O}(\mathbb{G}_a^r)\otimes A/A_{>0}\big)=\ker\big(A\to\mathcal{O}(U')\otimes A/A_{>0}\big)
				\end{align}
				i.e. the scheme theoretic images of $\mathbb{G}_a^r\times Z\to C$ and $U'\times Z\to C$ coincide. The $\subseteq$ direction is obvious. For $\supseteq$, let $h\in \ker\big(A\to \mathcal{O}(U')\otimes A/A_{>0}\big)$. By Proposition \ref{proposition: ideal of scheme theoretic image}, we have that $h\in A_{>0}$ and $\beta_1^{n_1}\cdots\beta_b^{n_b}.h\in A_{>0}$ for all $(n_1,\cdots,n_b)\in\mathbb{N}^b$. By Proposition \ref{proposition: ideal of scheme theoretic image} again, it suffices to prove $\gamma_1\cdots\gamma_m.h\in A_{>0}$ for all $m\in\mathbb{N}_+$ and all $\gamma_1,\cdots,\gamma_m\in\mathfrak{u}$. We have $\mathrm{Fit}_{d(\pi)-1}(\mathcal{Q}(\pi))\subseteq A_{>0}$ by Proposition \ref{proposition: dimStab and Fit}, and then the following is stronger 
				\begin{align}
					\gamma_1\cdots\gamma_m.h\in\langle\mathrm{U}(\mathfrak{v}).h\rangle+\mathrm{Fit}_{d(\pi)-1}(\mathcal{Q}(\pi))
				\end{align}
				where $\mathrm{U}(\mathfrak{v})$ is the universal enveloping algebra of $\mathfrak{v}$ and $\langle\mathrm{U}(\mathfrak{v}).h\rangle\subseteq A$ is the ideal generated by $\beta_1^{n_1}\cdots\beta_b^{n_b}.h$ for all $n\in\mathbb{N}^b$. Note that for all $g\in A$ and all $\zeta\in\mathfrak{u}$ we have $\zeta.g-\sum_{\mu=1}^b(\zeta.f_\mu)(\beta_\mu.g)\in\mathrm{Fit}_{d(\pi)-1}(\mathcal{Q}(\pi))$ (cf. Lemma \ref{lemma: Fit_(d-1) cap A_e is generated by Fit_(d-1) cap Fit_(r-1) cap A_e} (ii)). Then $\langle\mathrm{U}(\mathfrak{v}).h\rangle+\mathrm{Fit}_{d(\pi)-1}(\mathcal{Q}(\pi))\subseteq A$ is $\mathbb{G}_a^r$-invariant, and it contains $\zeta.h$ for all $\zeta\in\mathfrak{u}$. A simple induction on $m$ proves the statement. 
			\end{proof}
		}
	}

	\subsection{Statement and proof}
	{
		{
			\begin{proposition}\label{proposition: modification II}
				Let $\mathbb{G}_a^r\rtimes_w\mathbb{G}_m$ act on $\pi:C\to S$ as in Situation \ref{situation: G_a^r act on cone}. Assume 
				\begin{itemize}
					\item $\dim\mathrm{Stab}_{\mathfrak{u}}(z)=d(\pi)$ for all $z\in Z$; 
					\item $e(\pi)<\infty$, i.e. $\dim\mathrm{Stab}_{\mathfrak{u}}(x)<d(\pi)$ for some $x\in C$. 
				\end{itemize}
				Then the modification $\pi':C'\to S'$ with centre $\mathbb{G}_a^r.Z\hookrightarrow C$ satisfies 
				\begin{align}
					(d(\pi'),e(\pi'))<(d(\pi),e(\pi))\quad\textrm{in }\mathbb{N}\times\overline{\mathbb{N}_+}. 
				\end{align}
			\end{proposition}
			\begin{proof}
				Since $S$ is irreducible, the index $(d(\pi),e(\pi))$ does not change if we restrict to a non-empty open subvariety. We may assume conditions of Lemma \ref{lemma: Fit_(d-1) cap A_e is generated by Fit_(d-1) cap Fit_(r-1) cap A_e}: 
				\begin{itemize}
					\item $S$ is affine such that $\pi:C\to S$ corresponds to $A_0\to A=\bigoplus_{n\geq 0}A_n$; 
					\item there exist $\beta_1,\cdots,\beta_b\in\mathfrak{u}$ and $f_1,\cdots,f_b\in A_w$ for $b:=r-d(\pi)$ such that $\beta_i.f_j=\delta_{i,j}$. 
				\end{itemize}

				Let $J\subseteq A$ denote the ideal associated to the centre $\mathbb{G}_a^r.Z\hookrightarrow C$. Then $J$ is homogeneous and $\mathbb{G}_a^r$-invariant contained in $A_{>0}$. Moreover $J\ne0$ since $e(\pi)<\infty$. Let $i_1:=\min\{n:J_n\ne0\}\in\mathbb{N}_+$ be the minimal degree in $J$. In particular $J_{i_1}\subseteq A^{\mathfrak{u}}$. We have  
				\begin{align}
					S'=\mathrm{Proj}\big(A_0\oplus J_{i_1}\oplus J_{i_1}^2\oplus\cdots\big),\quad J_{i_1}^m:=\mathrm{Span}_{A_0}\{g_1\cdots g_m:g_1,\cdots,g_m\in J_{i_1}\}
				\end{align}
				which has an affine open covering 
				\begin{align}
					S'=\bigcup_{a\in J_{i_1}\setminus\{0\}}S'_a,\quad S'_a:=\mathrm{Spec}\Big(A_0\Big[\frac{J_{i_1}}{a}\Big]\Big). 
				\end{align}
				The affine cone $\pi':C'\to S'$ over $S'_a$ is 
				\begin{align}
					\pi'_a: C'_a\to S'_a,\quad C'_a:=\mathrm{Spec}\Big(A\Big[\frac{J}{a}\Big]\Big). 
				\end{align}

				It suffices to prove $(d(\pi'_a),e(\pi'_a))<(d(\pi),e(\pi))$. We can view $f_1,\cdots,f_b\in A\big[\frac{J}{a}\big]$ and then 
				\begin{align}
					1=\det(\beta_i.f_j)\in\mathrm{Fit}_{d(\pi)}(\mathcal{Q}(\pi'_a)). 
				\end{align}
				The integer $d(\pi'_a)\in\mathbb{N}$ is the smallest $d$ such that $\mathrm{Fit}_d(\mathcal{Q}(\pi_a'))\not\subseteq A\big[\frac{J}{a}\big]_{>0}$ by Proposition \ref{proposition: dimStab and Fit} (or \eqref{equation: Fit definition of d(pi)}). Then $d(\pi'_a)\leq d(\pi)$. If $d(\pi'_a)<d(\pi)$, then $(d(\pi'_a),e(\pi'_a))<(d(\pi),e(\pi))$. We then consider $d(\pi'_a)=d(\pi)$. 

				By Lemma \ref{lemma: Fit_(d-1) cap A_e is generated by Fit_(d-1) cap Fit_(r-1) cap A_e} (iii), there exists $h\in A_{e(\pi)+w}$ and $\zeta\in\mathfrak{u}$ such that 
				\begin{align}
					0\ne\zeta.h\in\mathrm{Fit}_{d(\pi)-1}(\mathcal{Q}(\pi))\cap A_{e(\pi)},\quad\beta_1.h=\cdots=\beta_b.h=0. 
				\end{align}
				We claim that $h\in J$. By Proposition \ref{proposition: ideal of scheme theoretic image}, it suffices to show $h\in A_{>0}$ and $\gamma_1\cdots\gamma_m.h\in A_{>0}$ for all $m\in\mathbb{N}_+$ and all $\gamma_1,\cdots,\gamma_m\in\mathfrak{u}$. By Lemma \ref{lemma: Fit_(d-1) cap A_e is generated by Fit_(d-1) cap Fit_(r-1) cap A_e} (i), we have $\gamma_m.h\in\mathrm{Fit}_{d(\pi)-1}(\mathcal{Q}(\pi))$. Since $\mathrm{Fit}_{d(\pi)-1}(\mathcal{Q}(\pi))$ is $\mathbb{G}_a^r$-invariant, we have $\gamma_1\cdots\gamma_m.h\in \mathrm{Fit}_{d(\pi)-1}(\mathcal{Q}(\pi))\subseteq A_{>0}$. Then $\frac{h}{a}\in A\big[\frac{J}{a}\big]$, and then $0\ne\frac{\zeta.h}{a}\in \mathrm{Fit}_{d(\pi'_a)-1}(\mathcal{Q}(\pi'_a))$, since \eqref{equation: nonzero in Fit_(d-1)} still holds. The minimal degree $e(\pi'_a)$ in $\mathrm{Fit}_{d(\pi'_a)-1}(\mathcal{Q}(\pi'_a))$ satisfies 
				\begin{align}
					e(\pi'_a)\leq\deg\Big(\frac{\zeta.h}{a}\Big)=e(\pi)-i_1<e(\pi). 
				\end{align}
			\end{proof}
		}
	}

	\subsection{Modifications in Example \ref{example: working example}}\label{subsection: working example - Mod II}
	{
		Assume $C^{(1)},C^{(3)}$ of Example \ref{example: working example} are known. We determine centres and describe $C^{(2)},C^{(4)}$. We apply previous calculations from Example \ref{example: working example - Q and Fit} and Example \ref{example: wokring example - index}. Note that we need $e(\pi)<\infty$ and $\dim\mathrm{Stab}_{\mathfrak{u}}(z)=d(\pi)$ for all $z\in Z$ to perform Modification II. These conditions are equivalent to the following 
		\begin{equation*}
			0\ne\mathrm{Fit}_{d(\pi)-1}(\mathcal{Q}(\pi))\subseteq\mathcal{A}_{>0},\quad \mathrm{Fit}_{d(\pi)}(\mathcal{Q}(\pi))=\mathcal{A}. 
		\end{equation*}
		These for $C^{(1)},C^{(3)}$ hold from Example \ref{example: working example - Q and Fit} and Example \ref{example: wokring example - index}. Recall the centre is $\mathbb{G}_a^2.Z$. 

		\subsubsection{$C^{(1)}$}
		{
			We have $\mathcal{A}^{(1)}=\Bbbk\big[a_{11},\tfrac{a_{21}}{a_{11}},\tfrac{a_{22}}{a_{11}},\tfrac{e}{a_{11}},\tfrac{f_1}{a_{11}},\tfrac{f_2}{a_{11}}\big]$. The zero section $Z^{(1)}$ is associated to the ideal 
			\begin{equation*}
				\mathcal{A}^{(1)}_{>0}=\big\langle\tfrac{a_{21}}{a_{11}},\tfrac{a_{22}}{a_{11}},\tfrac{e}{a_{11}},\tfrac{f_1}{a_{11}},\tfrac{f_2}{a_{11}}\big\rangle. 
			\end{equation*}
			The centre $\mathbb{G}_a^2.Z^{(1)}$ is associated to the ideal 
			\begin{equation*}
				\mathcal{J}^{(1)}=\ker\big(\mathcal{A}^{(1)}\to\mathcal{O}(\mathbb{G}_a^2)\otimes\mathcal{A}^{(1)}/\mathcal{A}^{(1)}_{>0}\big)=\big\langle\tfrac{a_{21}}{a_{11}},\tfrac{a_{22}}{a_{11}},\tfrac{e}{a_{11}},\tfrac{f_2}{a_{11}}\big\rangle. 
			\end{equation*}
			The homogeneous component of $\mathcal{J}^{(1)}$ of the minimal degree is $\mathcal{A}^{(1)}_0\tfrac{e}{a_{11}}$, with the minimal degree being $\sigma>0$. Then 
			\begin{equation*}\begin{split}
				\mathcal{A}^{(2)}={}&\mathcal{A}^{(1)}\Big[\frac{\mathcal{J}^{(1)}}{\tfrac{e}{a_{11}}}\Big]\\ 
				={}&\textstyle{\Bbbk\big[a_{11},\tfrac{a_{21}}{e},\tfrac{a_{22}}{e},\tfrac{e}{a_{11}},\tfrac{f_1}{a_{11}},\tfrac{f_2}{e}\big]}
			\end{split}\end{equation*}
		}

		\subsubsection{$C^{(3)}$}
		{
			We have $\mathcal{A}^{(3)}=\Bbbk\big[a_{11},\tfrac{a_{21}}{e},\tfrac{a_{22}}{e},\tfrac{e}{a_{11}},\tfrac{f_1}{a_{11}},\tfrac{f_2}{e}\big]$. The zero section $Z^{(3)}$ is associated to the ideal 
			\begin{equation*}
				\mathcal{A}^{(3)}_{>0}=\big\langle\tfrac{a_{21}}{e},\tfrac{a_{22}}{e},\tfrac{e}{a_{11}},\tfrac{f_1}{a_{11}},\tfrac{f_2}{e}\big\rangle. 
			\end{equation*}
			The centre $\mathbb{G}_a^2.Z^{(3)}$ is associated to the ideal 
			\begin{equation*}
				\mathcal{J}^{(3)}=\ker\big(\mathcal{A}^{(3)}\to\mathcal{O}(\mathbb{G}_a^2)\otimes\mathcal{A}^{(3)}/\mathcal{A}^{(3)}_{>0}\big)=\big\langle\tfrac{a_{21}}{e},\tfrac{a_{22}}{e},\tfrac{e}{a_{11}},\tfrac{f_2}{e}\rangle. 
			\end{equation*}
			The homogeneous component of $\mathcal{J}^{(3)}$ of the minimal degree is $\mathcal{A}^{(3)}_0\tfrac{a_{21}}{e}+\mathcal{A}^{(3)}_0\tfrac{a_{22}}{e}$, with the minimal degree being $\rho-\sigma$. Then $C^{(4)}=\mathrm{Spec}(\mathcal{A}^{(4),1})\bigcup\mathrm{Spec}(\mathcal{A}^{(4),2})$ is covered by two affine opens with 
			\begin{equation*}\begin{split}
				\mathcal{A}^{(4),1}={}&\mathcal{A}^{(3)}\Big[\frac{\mathcal{J}^{(3)}}{\tfrac{a_{21}}{e}}\Big]\\ 
				={}&\textstyle{\Bbbk\big[a_{11},\tfrac{a_{21}}{e},\tfrac{a_{22}}{a_{21}},\tfrac{e^2}{a_{11}a_{21}},\tfrac{f_1}{a_{11}},\tfrac{f_2}{a_{21}}\big]}\\ 
				\mathcal{A}^{(4),2}={}&\mathcal{A}^{(3)}\Big[\frac{\mathcal{J}^{(3)}}{\tfrac{a_{22}}{e}}\Big]\\ 
				={}&\textstyle{\Bbbk\big[a_{11},\tfrac{a_{21}}{a_{22}},\tfrac{a_{22}}{e},\tfrac{e^2}{a_{11}a_{22}},\tfrac{f_1}{a_{11}},\tfrac{f_2}{a_{22}}\big]}. 
			\end{split}\end{equation*}

		}
	}
}

\appendix
\section{Ideal description}
{
	\subsection{Secondary Fitting ideal}
	{
		Recall $e(\pi)$ defined as the minimal degree in $\mathrm{Fit}_{d(\pi)-1}(\mathcal{Q}(\pi))$ in Definition \ref{definition: index (d,e)}. The following lemma describes $\mathrm{Fit}_{d(\pi)-1}(\mathcal{Q}(\pi))\cap\mathcal{A}_{e(\pi)}$ locally, tracing $e(\pi)$ along modifications. 
		{
			\begin{lemma}\label{lemma: Fit_(d-1) cap A_e is generated by Fit_(d-1) cap Fit_(r-1) cap A_e}
				Let $\mathbb{G}_a^r\times_w\mathbb{G}_m$ act on $\pi:C\to S$ as in Situation \ref{situation: G_a^r act on cone}. Assume 
				\begin{itemize}
					\item $\dim\mathrm{Stab}_{\mathfrak{u}}(z)=d(\pi)$ for all $z\in Z$; 
					\item $e(\pi)<\infty$, i.e. $\dim\mathrm{Stab}_{\mathfrak{u}}(x)<d(\pi)$ for some $x\in C$; 
					\item $S$ is affine such that $\pi:C\to S$ corresponds to $A_0\to A=\bigoplus_{n\geq0}A_n$; 
					\item there exist $\beta_1,\cdots,\beta_b\in\mathfrak{u}$ and $f_1,\cdots,f_b\in A_w$ for $b:=r-d(\pi)$ such that $\beta_i.f_j=\delta_{i,j}$ for $1\leq i,j\leq b$. 
				\end{itemize}
				Then 
				\begin{itemize}
					\item[(i)] For $\zeta\in\mathfrak{u}$ and $h\in A$, we have $\zeta.h-\sum_{\mu=1}^b(\zeta.f_\mu)(\beta_\mu.h)\in\mathrm{Fit}_{d(\pi)-1}(\mathcal{Q}(\pi))$; 
					\item[(ii)] For $h\in A$, we have $\beta_1.\hat{h}=\cdots=\beta_b.\hat{h}=0$ for
					\begin{align}
						\hat{h}:=\sum_{n\in\mathbb{N}^b}\frac{(-1)^{|n|}}{n!}(\beta^n.h)f^n=h+\sum_{n\in\mathbb{N}^b\setminus\{0\}}\frac{(-1)^{|n|}}{n!}(\beta^n.h)f^n
					\end{align}
					where $|n|:=n_1+\cdots+n_b$ and $n!:=n_1!\cdots n_b!$ and $\beta^n:=\beta_1^{n_1}\cdots\beta_b^{n_b}$ and $f^n:=f_1^{n_1}\cdots f_b^{n_b}$; 
					\item[(iii)] The $A_0$-module $\mathrm{Fit}_{d(\pi)-1}(\mathcal{Q}(\pi))\cap A_{e(\pi)}$ is generated by the set of elements 
					\begin{align}
						\left\{\zeta.h:\begin{matrix}\zeta\in\mathfrak{u},\textrm{ and }h\in A_{e(\pi)+w}\\\textrm{such that }\beta_1.h=\cdots=\beta_b.h=0\end{matrix}\right\}. 
					\end{align}
				\end{itemize}
			\end{lemma}
			\begin{proof}
				Denote $\mathfrak{v}:=\mathrm{Span}_\Bbbk\{\beta_1,\cdots,\beta_b\}\subseteq \mathfrak{u}$ and $A^{\mathfrak{v}}:=\{h\in A:\beta_1.h=\cdots=\beta_b.h=0\}$. Then (i) is from 
				\begin{align}
					\zeta.h-\sum_{\mu=1}^b(\zeta.f_\mu)(\beta_\mu.h)=\det\begin{pmatrix}
					\beta_1.f_1&\cdots&\beta_1.f_b&\beta_1.h\\
					\vdots&\ddots&\vdots&\vdots\\
					\beta_b.f_1&\cdots&\beta_b.f_b&\beta_b.h\\
					\zeta.f_1&\cdots&\zeta.f_b&\zeta.h
					\end{pmatrix}\in\mathrm{Fit}_{d(\pi)-1}(\mathcal{Q}(\pi)); 
				\end{align}
				and (ii) is by direct calculation. 

				The $A_0$-module $\mathrm{Fit}_{d(\pi)-1}(\mathcal{Q}(\pi))\cap A_{e(\pi)}$ is generated by elements of the form 
				\begin{align}
					\det\begin{pmatrix}
						\zeta_1.h_1&\cdots&\zeta_1.h_{b+1}\\
						\vdots&\ddots&\vdots\\
						\zeta_{b+1}.h_1&\cdots&\zeta_{b+1}.h_{b+1}
					\end{pmatrix}\ne0
				\end{align}
				for $\zeta_1,\cdots,\zeta_{b+1}\in\mathfrak{u}$ and $h_1,\cdots,h_{b+1}\in\mathcal{A}_{\geq w}$ homogeneous such that $\det(\zeta_i.h_j)\in A_{e(\pi)}$. We can assume that $h_1,\cdots,h_{b+1}$ have non-decreasing degrees. It suffices to prove $\det(\zeta_i.h_j)$ can be generated by the set of (iii). We have $h_{b+1}\in\mathcal{A}_{>w}$, otherwise $\det(\zeta_i.h_j)\in \mathrm{Fit}_{d(\pi)-1}(\mathcal{Q}(\pi))\cap A_0=0$. If $h_j\in\mathcal{A}_{>w}$ and $n\in\mathbb{N}^b\setminus\{0\}$ with $n_\mu>0$, then 
				\begin{equation}\begin{split}\label{equation: det(zeta h)=det(zeta hat h)}
					&\quad\det\begin{pmatrix}
						\zeta_1.h_1&\cdots&\zeta_1.((\beta^n.h_j)f^n)&\cdots&\zeta_1.h_{b+1}\\
						\vdots&\vdots&\vdots&\vdots&\vdots\\
						\zeta_{b+1}.h_1&\cdots&\zeta_{b+1}.((\beta^n.h_j)f^n)&\cdots&\zeta_{b+1}.h_{b+1}
					\end{pmatrix}\\
					&= f_\mu\begin{pmatrix}
						\zeta_1.h_1&\cdots&\zeta_1.((\beta^n.h_j)f^{n-\epsilon_\mu})&\cdots&\zeta_1.h_{b+1}\\
						\vdots&\vdots&\vdots&\vdots&\vdots\\
						\zeta_{b+1}.h_1&\cdots&\zeta_{b+1}.((\beta^n.h_j)f^{n-\epsilon_\mu})&\cdots&\zeta_{b+1}.h_{b+1}
					\end{pmatrix}\\
					&\quad+((\beta^n.h_j)f^{n-\epsilon_\mu})\begin{pmatrix}
						\zeta_1.h_1&\cdots&\zeta_1.f_\mu&\cdots&\zeta_1.h_{b+1}\\
						\vdots&\vdots&\vdots&\vdots&\vdots\\
						\zeta_{b+1}.h_1&\cdots&\zeta_{b+1}.f_\mu&\cdots&\zeta_{b+1}.h_{b+1}
					\end{pmatrix}\\
					&= 0+0
				\end{split}\end{equation}
				where $\epsilon_\mu=(0,\cdots,1,\cdots,0)\in\mathbb{N}^b$ is the $\mu$th basis vector, and the two determinants vanish since they are in $\mathrm{Fit}_{d(\pi)-1}(\mathcal{Q}(\pi))\cap A_{e(\pi)-w}$ and $\mathrm{Fit}_{d(\pi)-1}(\mathcal{Q}(\pi))\cap A_{e(\pi)-(\deg h_j-w)}$ respectively. 

				Define 
				\begin{align}
					\hat{h}_j:=\begin{cases}
						h_j,&\textrm{if }h_j\in A_w\\
						h_j+\sum_{n\in\mathbb{N}^b\setminus\{0\}}\frac{(-1)^{|n|}}{n!}(\beta^n.h_j)f^n,&\textrm{if }h_j\in A_{>w}. 
					\end{cases}
				\end{align}
				We have 
				\begin{itemize}
					\item if $\hat{h}_j\in A^{\mathfrak{v}}$, then $\zeta_i.\hat{h}_j\in\mathrm{Fit}_{d(\pi)-1}(\mathcal{Q}(\pi))$ by (i); 
					\item $\hat{h}_j\in A_w$ or $\hat{h}_j\in A^{\mathfrak{v}}\cap A_{>w}$ by (ii); 
					\item $\det(\zeta_i.h_j)=\det(\zeta_i.\hat{h}_j)$ by \eqref{equation: det(zeta h)=det(zeta hat h)}. 
				\end{itemize}
				The degree of $\det(\zeta_i.h_j)=\det(\zeta_i.\hat{h}_j)$ is 
				\begin{equation}\begin{split}
					e(\pi)=\sum_{j=1}^{b+1}\deg(\zeta_i.\hat{h}_j)=\sum_{j:\deg(h_j)>w}\deg(\zeta_i.\hat{h}_j)\geq\sum_{j:\deg(h_j)>w}e(\pi)
				\end{split}\end{equation}
				where the inequality is from $\zeta_i.\hat{h}_j\in \mathrm{Fit}_{d(\pi)-1}(\mathcal{Q}(\pi))\subseteq A_{\geq e(\pi)}$. Then $h_1,\cdots,h_b\in A_w$ and $h_{b+1}\in A_{e(\pi)+w}$. Then $\det(\zeta_i.h_j)$ is in the $A_0$-span of the set of (iii), since $\det(\zeta_i.h_j)=\det(\zeta_i.\hat{h}_j)=\sum_{i=1}^{b+1}a_i(\zeta_i.\hat{h}_{b+1})$, where $a_i\in A_0$ are cofactors and $\hat{h}_{b+1}\in A_{e(\pi)+w}\cap A^{\mathfrak{v}}$. 
			\end{proof}
		}
	}
	\subsection{Scheme theoretic image}
	{
		Let $\mathbb{G}_a^r$ act on $\mathrm{Spec}(R)$ via the co-action map $\sigma:R\to R[u_1,\cdots,u_r]$, where $u_1,\cdots,u_r\in\mathfrak{u}^*$ form a basis. Let $\xi_1,\cdots,\xi_r\in\mathfrak{u}$ be the dual basis. Let $I\subseteq R$ be an ideal. Then the scheme theoretic image of $\mathbb{G}_a^r\times \mathbb{V}(I)\to \mathrm{Spec}(R)$ is the subscheme associated to the ideal $\ker(R\to R/I[u_1,\cdots,u_r])$. 
		{
			\begin{proposition}\label{proposition: ideal of scheme theoretic image}
				Let $\mathbb{G}_a^r\curvearrowright\mathrm{Spec}(R)$, $I\subseteq R$ be as above. For $h\in R$, the following are equivalent: 
				\begin{itemize}
					\item $h\in \ker(R\to R/I[u_1,\cdots,u_r])$; 
					\item $\xi^n.h\in I$ for all $n\in\mathbb{N}^r$, where $\xi^n:=\xi_1^{n_1}\cdots\xi_r^{n_r}$. 
				\end{itemize}
			\end{proposition}
			\begin{proof}
				This is immediate from the Taylor expansion 
				\begin{align}
					\sigma(h)=\sum_{n\in\mathbb{N}^r}\frac{\xi^n.h}{n!}u^n,\quad h\in R
				\end{align}
				where $n!:=n_1!\cdots n_r!$ and $u^n:=u_1^{n_1}\cdots u_r^{n_r}$. 
			\end{proof}
		}
	}

}

\printbibliography
\end{document}